\documentclass[onefignum,onetabnum]{siamart190516}

\usepackage{subcaption}
\usepackage{lipsum}
\usepackage{amsfonts}
\usepackage{graphicx}
\usepackage{epstopdf}
\usepackage{algorithmic}
\ifpdf
\DeclareGraphicsExtensions{.eps,.pdf,.png,.jpg}
\else
\DeclareGraphicsExtensions{.eps}
\fi

\newsiamremark{remark}{Remark}
\newsiamremark{hypothesis}{Hypothesis}
\crefname{hypothesis}{Hypothesis}{Hypotheses}
\newsiamthm{claim}{Claim}

\headers{Fast variable density 3-D node generation}{K. van der Sande, and B. Fornberg}

\title{Fast variable density 3-D node generation\thanks{Submitted to the editors May 9, 2020.
		\funding{This research did not receive any specific grant from funding agencies in the public, commercial, or not-for-profit sectors.}}}

\author{Kiera van der Sande\thanks{Department of Applied Mathematics, University of Colorado Boulder, Boulder, CO
		(\email{kiera.vandersande@colorado.edu}).}
	\and Bengt Fornberg\thanks{Department of Applied Mathematics, University of Colorado Boulder, Boulder, CO 
		(\email{fornberg@colorado.edu}).}
}

\usepackage{amsopn}

\ifpdf
\hypersetup{
	pdftitle={Fast variable density 3-D node generation},
	pdfauthor={Kiera van der Sande, and Bengt Fornberg}
}
\fi

\begin{document}

\maketitle
	
\begin{abstract}
	Mesh-free solvers for partial differential equations perform best on scattered quasi-uniform nodes. Computational efficiency can be improved by using nodes with greater spacing in regions of less activity. However, there is no ideal way to generate nodes for these solvers. We present an advancing front type method to generate variable density nodes in 2-D and 3-D with clear generalization to higher dimensions. The exhibited cost of generating a node set of size $N$ in 2-D and 3-D with the present method is $O(N)$.            
\end{abstract}

\begin{keywords}
	Node generation,  variable density points, mesh-free PDE solvers, RBF-FD
\end{keywords}

\begin{AMS}
	65N99, 65Y20, 65M70
\end{AMS}

\section{Introduction}

Mesh-free methods for solving partial differential equations such as radial basis function-generated finite differences (RBF-FD) have become increasingly popular. These methods use scattered nodes of variable density rather than a mesh as a computational domain. RBF-FD methods allow for high geometric flexibility, but still require certain constraints on node sets in order to ensure solution accuracy and stability \cite{FF2015primer,FF2015actanumerica}. For example, nodes that are locally too irregular can be problematic for the stability of PDE solvers. Hence, one key quality requirement is for nodes to be locally quasi-uniform i.e., if you zoom into a region of nodes they should be close to equispaced (see \cref{sec:nodesetquality} for more rigorous definitions of node quality). Node generation algorithms should also satisfy minimum spacing and bounded gaps between nodes and their neighbors, and have the ability to spatially vary node density in a prescribed manner. In general, placing nodes in a domain is reminiscent of circle packing in 2-D and sphere packing in 3-D. Optimal node sets for a given node spacing should be as densely packed as possible while maintaining a prescribed distance between nodes. 

Node generation remains a bottleneck for mesh-free PDE solvers, especially in higher than 2 dimensions or where variable density is desired. Recent work has been done on producing quality node sets specifically for RBF-FD \cite{Flyer2010,Fornberg2015nodegen,Slak2019,Vlasiuk2018}. Here we build off the method of 2-D node generation from Fornberg \& Flyer \cite{Fornberg2015nodegen} to generate nodes in higher dimensions. This previous method was constrained to 2 dimensions due to the way that nodes were generated and stored. The present method utilizes a background grid and local searches to allow quality nodes of variable density to be generated in 2-D or 3-D according to a desired node spacing function, with the ability to generalize to higher dimensions. The algorithm also guarantees a minimum spacing requirement between nodes. In this work we do not seek to further demonstrate the robustness of RBF-FD and other mesh-free PDE methods, rather to fill the need for locally quasi-uniform and variable density node sets that these methods require \cite{Ahmad2020local,Flyer2010,FF2015primer}. 

Current methods of node generation can be categorized broadly into iterative methods, sphere packing methods and advancing front methods. Often unstructured grid generators are used and nodes extracted as the vertices of the mesh. This is, however, computationally wasteful since RBF-FD methods make no use of the often costly step of connecting nodes into good aspect ratio elements. Iterative methods begin with an initial node set and update their positions through either a form of energy minimization \cite{hardin2016comparison}, short-range interaction forces \cite{nie2014parallel} or gradient flow \cite{Vlasiuk2018}. These methods are strongly dependent on their initial configuration and can be costly. Sphere and circle packing methods can be extended to arbitrary dimension and can be parallelizable, but are often constrained to constant sphere radii (i.e. constant node density). These methods include Poisson disk sampling \cite{bridson2007fast,dunbar2006spatial,talmor1997well}, which is used for sampling in the graphics community and was recently introduced as a method of generating nodes for RBF-FD \cite{Shankar2018}. There has been some interest in defining the requirements for variable radii in this context \cite{mitchell2012variable}. Advancing front methods are computationally more efficient and relatively simple to implement. The proposed method in this paper is an advancing front type method, as is the original in \cite{Fornberg2015nodegen} and the work of \cite{li2000biting,Lohner1998,Slak2019}. 

The rest of the paper is organized as follows: \cref{sec:algorithm} outlines the method in 3-D and presents two possible modifications to the basic method; \cref{sec:generatednodesets} investigates different metrics of node quality and compares the present method to other node sets in 2-D and 3-D; and \cref{sec:RBFmethods} demonstrates the application of the method for use in RBF-FD. Conclusions and future work are presented in \cref{sec:Conclusion}.

\section{Node generation in arbitrary dimension}
\label{sec:algorithm}

\subsection{The basic node generation algorithm}
We outline the algorithm for generating node sets by first considering the 3-D case in a bounded box. The desired spatial density of the nodes is specified through an exclusion radius function $r(x,y,z)$, which can be any 3-D function and defines the minimum spacing between nodes.

The method is an advancing-front type method, which relies on a background grid. In 3-D this is a dense grid in one Cartesian plane and the front progresses in the normal direction to it. For simplicity the grid is considered to be in the $x,y$-plane and the front to move in the increasing $z$-direction. The grid is stored as an array of `potential dot placements' (PDPs) with associated `heights' in the normal direction. These heights are initialized as the bottom $z$-plane of the boundary box plus a small random perturbation (on the order of the minimum desired distance between nodes at $z=0$). The first placed node is chosen as the minimum of these heights. The method proceeds as described in \cref{alg:method}. 

\vspace{0.5em}
\begin{algorithm}
	\caption{Node generation pseudocode}
	\label{alg:method}
\begin{algorithmic}
	\STATE{Initialize PDP array to the height of the bottom of the bounding box.}
	\STATE{Choose the minimum of the initial array as the first node location $p$.}
	\WHILE{the lowest PDP height is within the bounding box}
	\STATE{Add $p$ to the list  of generated nodes.}
	\STATE{Calculate the exclusion radius $r(p)$ at $p$.}
	\IF{updated height $>$ current height}
	\STATE{Update the heights of PDP within the sphere of radius $r$ centered at $p$ to lie on the upper half of this sphere.}
	\ENDIF
	\STATE{Set $x_0$ to be the PDP location with the minimum updated height.}
	\STATE{Set $x_1$ to be PDP location with the minimum height within $2r(p)$ of $x_0$.} 
	\WHILE{$|x_{i+1}-x_{i}| > r(p)$}
		\STATE{Update $x_i = x_{i+1}$}
		\STATE{Set $x_{i+1}$ to be the PDP location with minimum height within $2r(p)$ of $x_i$.}
	\ENDWHILE
	\STATE{Let the next node location be $p = x_{i+1}$.}
	\ENDWHILE
\end{algorithmic}
\end{algorithm}	
\vspace{0.5em}

\begin{figure}[htbp]
	\centering
	\begin{subfigure}[t]{0.48\textwidth}
		\centering
		\includegraphics[width=0.9\textwidth]{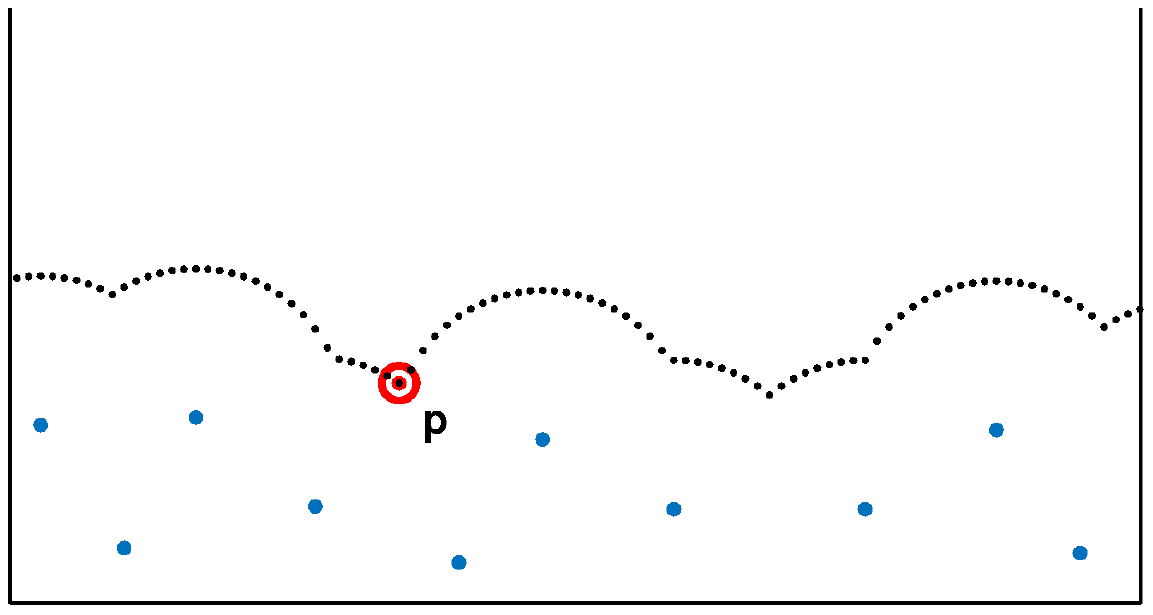}
		\caption{$p$ is the current node}
	\end{subfigure}
	\begin{subfigure}[t]{0.48\textwidth}
		\centering
		\includegraphics[width=0.9\textwidth]{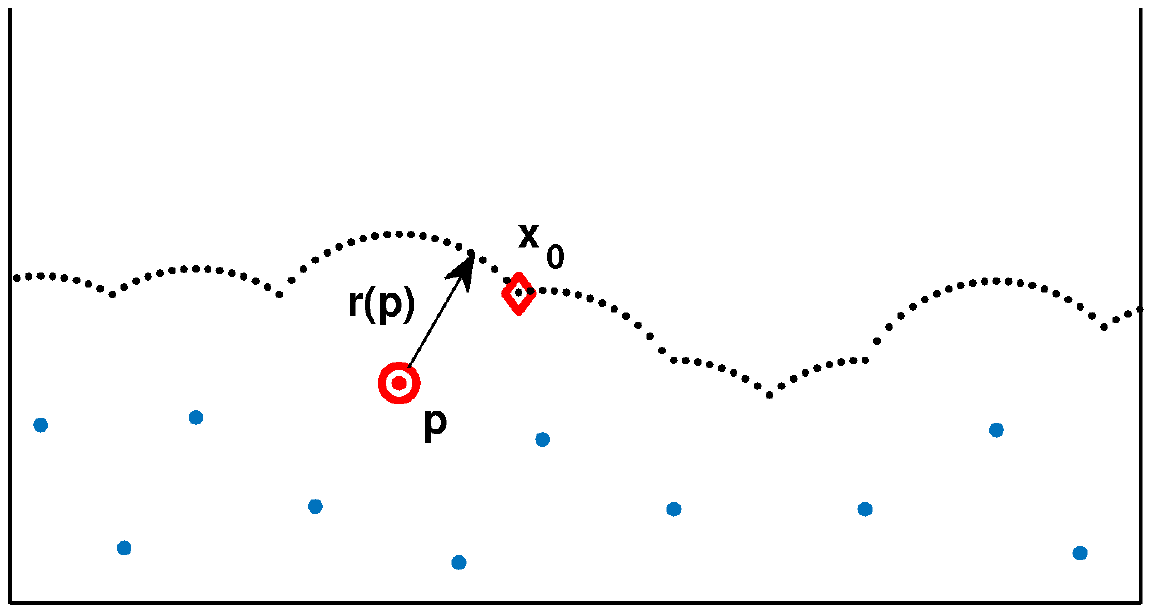}
		\caption{PDP heights are updated}
	\end{subfigure}
	\begin{subfigure}[t]{0.48\textwidth}
		\centering
		\includegraphics[width=0.9\textwidth]{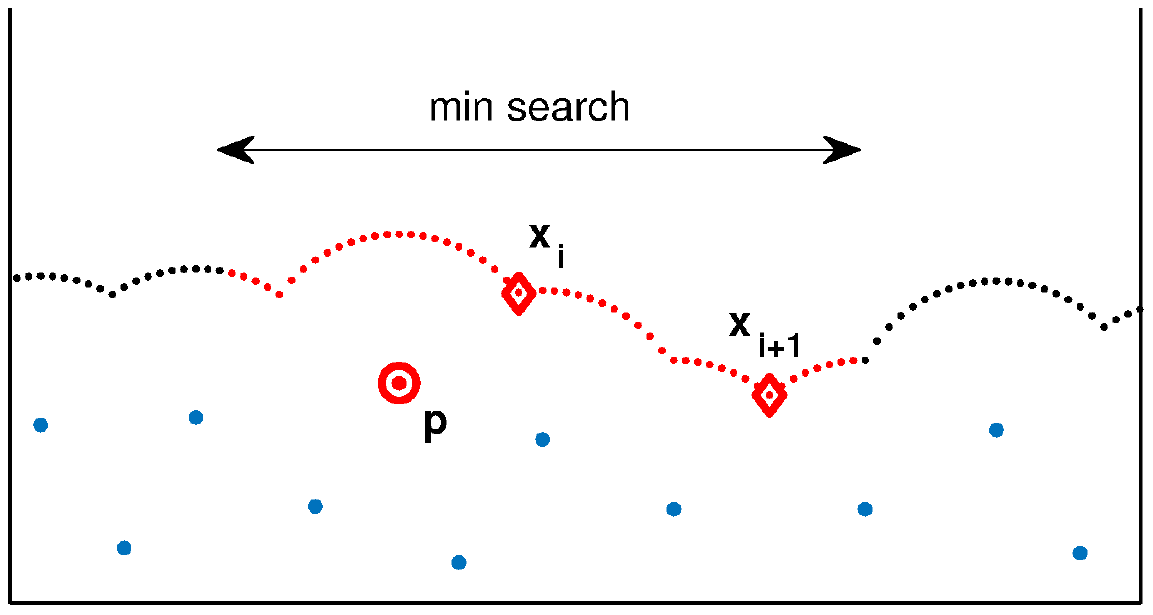}
		\caption{Search for local minimum}
	\end{subfigure}
	\begin{subfigure}[t]{0.48\textwidth}
		\centering
		\includegraphics[width=0.9\textwidth]{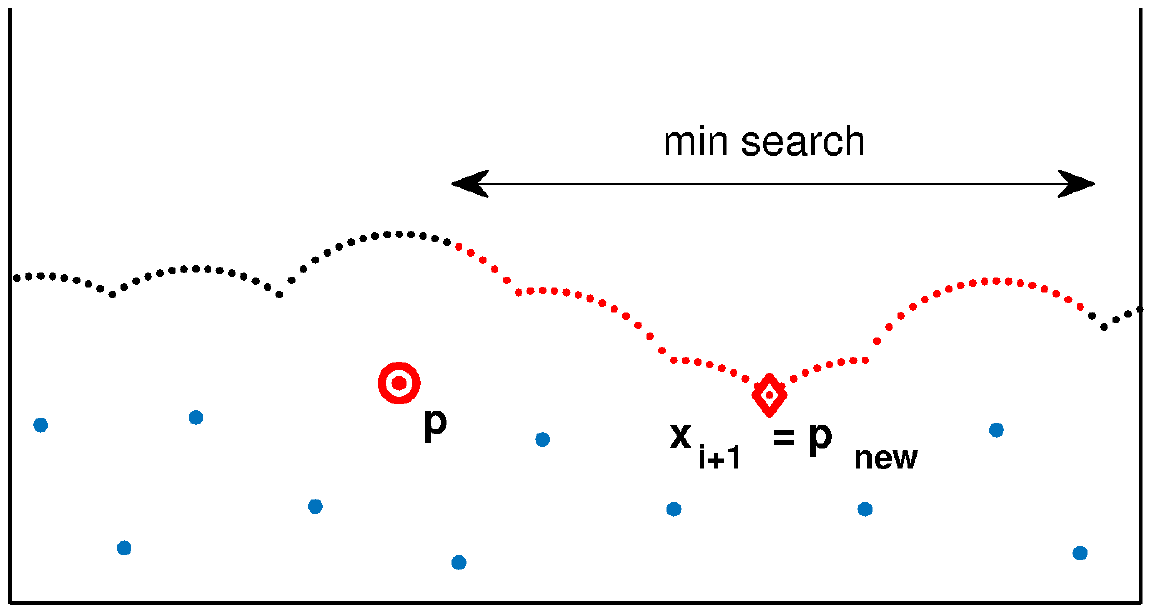}
		\caption{$p_{new}$ is the next node.}
	\end{subfigure}
	\caption{Illustration of node generation algorithm in 2-D. The `potential dot placements' (PDPs) shown by small black dots are are an advancing front. }
	\label{fig:advancingfront}
\end{figure}

In 3-D and higher there is no way to sort the PDPs to track the global minimum, as there was in the 2-D method of \cite{Fornberg2015nodegen}.  In order to find a close local minimum to the last placed node $p$, an iterative moving window search over the PDPs is used. The minimum of the updated heights is set as $x_0$ and iterations are taken from there to find a local minimum. At each iteration, $x_{n+1}$ is set to be the minimum of the PDPs within a radius of $2r(p)$ of $x_n$. This continues until $x_{n+1}$ is within $r(p)$ of $x_n$. If the edge of the box is reached, the search wraps around to the other side of the domain. \texttt{MATLAB} code for the algorithm is provided in \cite{nodegencode}. A visual representation of the algorithm in 2-D is shown in \cref{fig:advancingfront}.

Note that the resolution of the background grid will have an effect on the resulting node set. The grid should be fine enough to resolve spatial varying node densities, but not so fine as to impede efficiency. Experiments showed that setting $\Delta x$ of the background grid to be 10 times smaller than the minimum desired exclusion radius is a good compromise. Lower grid density gave way to discretization errors in finding local minima so nodes ended up further apart than desired, while increasing grid density above a factor of around $10-20$ did not significantly change the resulting node quality. Timing tests also showed that reducing the background grid density further gave no additional benefit as the cost is dominated by other factors. Unless otherwise stated, a factor of 10 will be used throughout this work.

We are interested in interior nodes so we generate nodes in a box (or sphere, as described later in \cref{sec:Spherical}). In order to create node sets in more complex domains, nodes should first be generated in a bounding box and points outside the desired domain discarded after the box is filled. If desired, boundary nodes can then be added manually. A final step would simulate local electrostatic repulsion on nodes near the boundary, as described in \cite{Fornberg2015nodegen} and utilized in \cite{Ahmad2020local,Flyer2010}. Treatment of boundaries in RBF-FD has been described previously in \cite{Bayona2019role,Bayona2017role}.

Generating nodes in higher dimension $d$ will require a PDP array of dimension $d-1$ and the front will advance in the last dimension. For example in 4-D, the background PDP grid will be a 3-D array and nodes will be placed in increasing `height' in the 4th dimension. 
    
\subsection{Spatial density function and exclusion radius}
The exclusion radius is given by a 3-D function that prescribes the desired distance between nodes. Setting a uniform exclusion radius $r$ for this method will guarantee a minimum spacing of $r$ between nodes. In the case of a spatially varying exclusion radius, a minimum spacing requirement can be defined based on the exclusion radius of the previously placed nodes (\textit{prior-disks}), the current node (\textit{current-disks}), or some function of the two (See \cref{sec:nodesetquality} for further details). 

The present method naturally adheres to the prior-disks variation as the exclusion radii of the previously placed nodes defines the position that the current node will be placed at. We propose a couple ways to deal with variations of this minimal spacing requirement in \cref{sec:directioncorrection} and \cref{sec:Spherical}.

\subsection{Direction dependence and a possible correction}
\label{sec:directioncorrection}

In the algorithm each node is placed based on the exclusion radii of the previously placed nodes and then the front is updated to include the exclusion radii of the new node. However this may lead to a directional dependence since the front advances in one particular direction (i.e. the Cartesian $z$-direction). An exclusion radius function which varies in $z$ will then have some systematic error in the $z$-direction. There are several possible ways to correct for this direction dependence and satisfy different minimal spacing requirements. Here we introduce one possible correction in the spirit of a \textit{bigger-disks} minimal spacing. 

A first order correction is added when placing a new node $p_j$. Before placing the node, a check is performed of whether any already placed nodes are within the exclusion radius $r(p_j)$ of the new node. In order to avoid an expensive search of the list of previously generated nodes we store an additional array, which is the same size as the background grid and initialized with null values. Each element in this array corresponds to an $(x,y)$ location in the background grid. When a node is placed at a given $(x,y)$ location a pointer to that node's place in the list of generated nodes is stored in the corresponding element of this additional array. This allows for a check of only close enough background grid elements to see if there are any placed nodes nearby. If any nodes are within $r(p_j)$, the height of the new node is increased until its own exclusion radius is satisfied. In 3-D, this correction can be written as:
\begin{equation}
z_{new} = z_{nbr} + \sqrt{r^2 - (x-x_{nbr})^2-(y-y_{nbr})^2}
\label{eqn:zcorrection}
\end{equation}
where $(x_{nbr},y_{nbr},z_{nbr})$ is the position of the nearest neighbor. The new node is then placed at $(x,y,z_{new})$.

\subsection{A modification for spherical density functions}
\label{sec:Spherical}
Often it is desirable to have node density vary in the radial direction, i.e. in modeling an atom or the atmosphere. Another option for avoiding direction dependence in this case is to construct the node set as a radially advancing front. This method allows better control of whether the minimal spacing requirement should be based on the minimum, maximum, or average exclusion radius between the current node to be placed and the previously placed nodes. When placing the next node based on the exclusion radii of previous nodes, the bias will be in the radial direction rather than the $z$-direction, which will allow for better radial symmetry.

To implement this modification, instead of starting with a background grid in the $x,y$-plane a grid can be constructed in $(\theta,\phi)$ and the algorithm can be carried out in spherical coordinates building outwards from the origin. 

\section{Generated node sets}
\label{sec:generatednodesets}
Nodes are generated in 2-D and 3-D, and compared to existing algorithms. For the measures of node quality used in the following section, we refer to \cref{sec:nodesetquality}.

\subsection{Nodes in 2-D}
\label{sec:2Dresults}
\begin{figure}[htbp]
	\centering
	\begin{subfigure}{0.32\textwidth}
		\centering
		\includegraphics[width=0.9\textwidth]{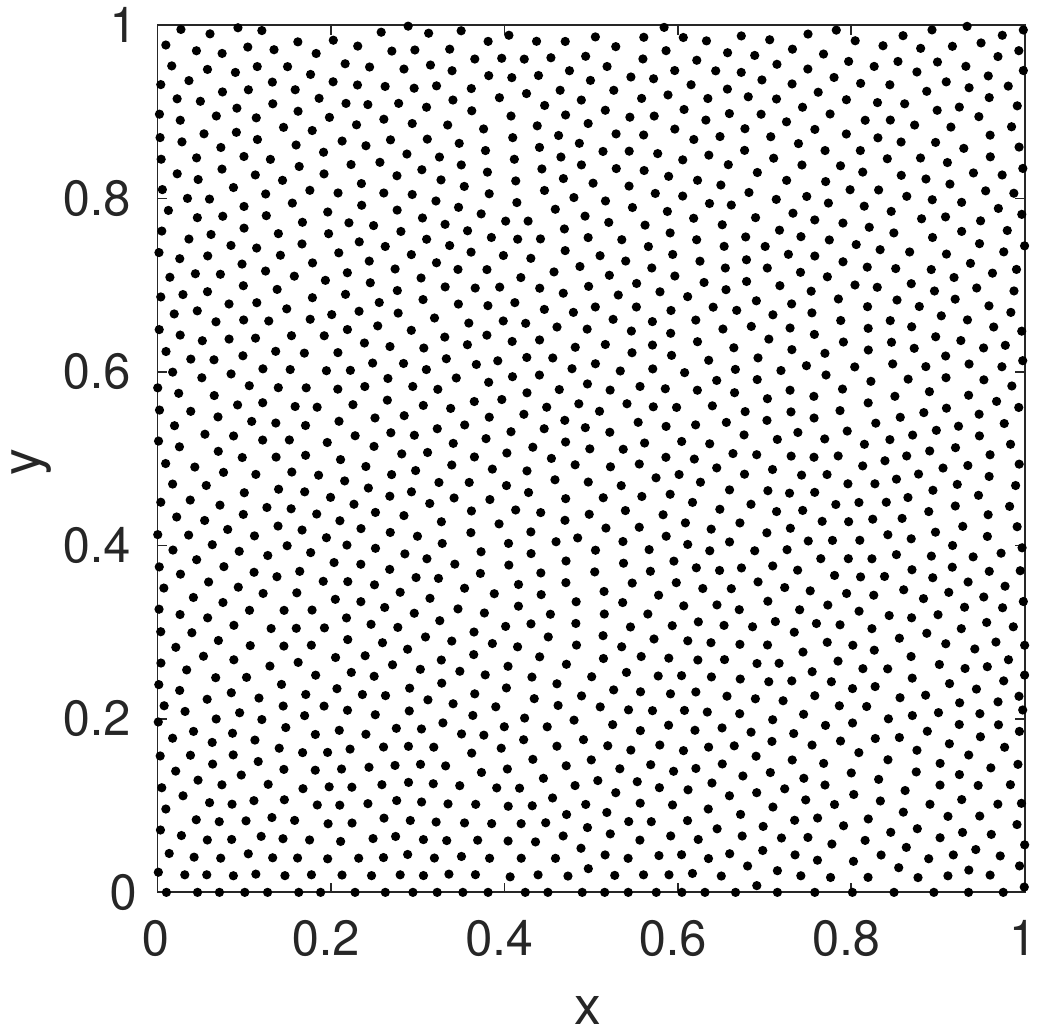}
		\caption{\centering
			Fornberg \& Flyer, $N = 1492$}
		\label{fig:FFnodes}
	\end{subfigure}
	\begin{subfigure}{0.32\textwidth}
		\centering
		\includegraphics[width=0.9\textwidth]{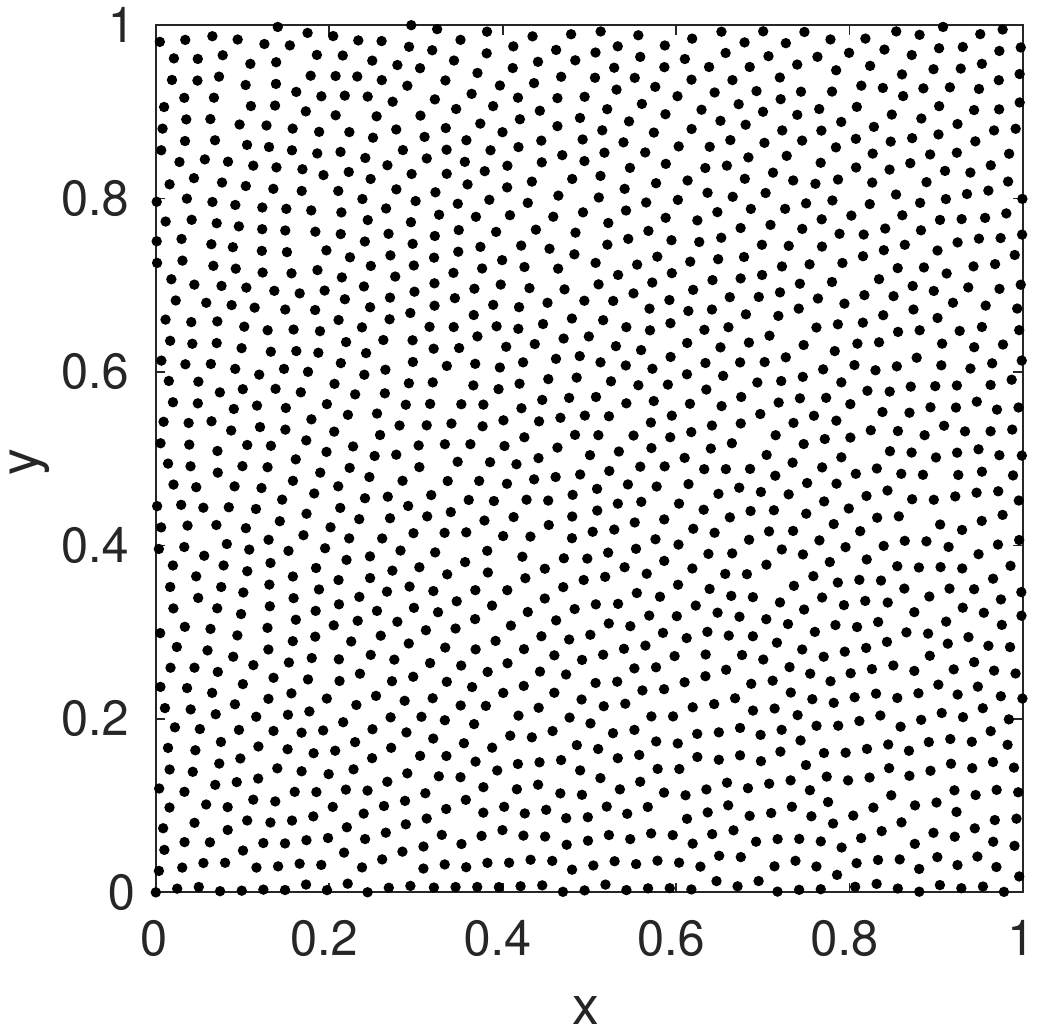}
		\caption{Slak \& Kosec,\\ $N = 1429$}
		\label{fig:SKnodes}
	\end{subfigure}
	\begin{subfigure}{0.32\textwidth}
		\centering
		\includegraphics[width=0.9\textwidth]{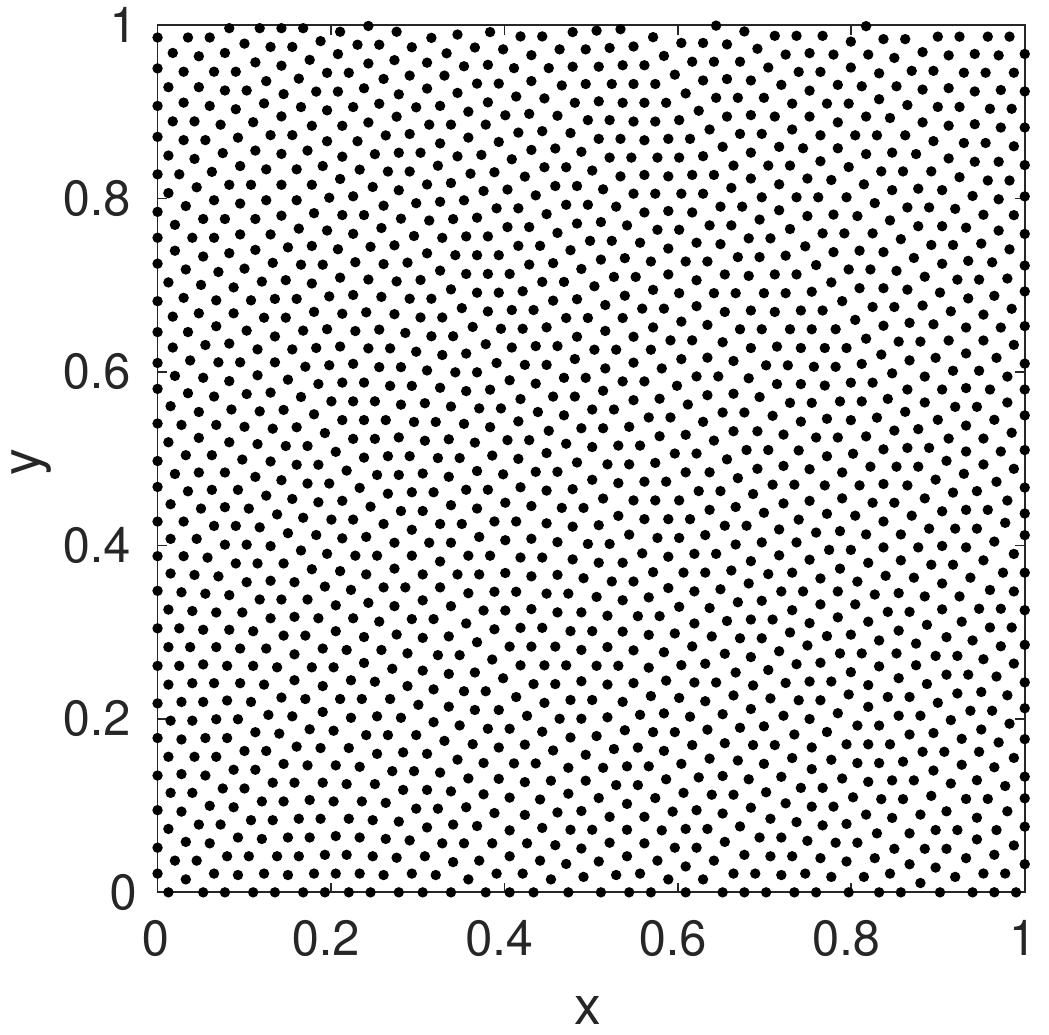}
		\caption{\centering Present Method, $N = 1698$}
		\label{fig:newnodes}
	\end{subfigure}
	\caption{2-D uniform node sets with $r=0.025$ spacing, $n$ being the the number of sample points created at each step of the algorithms in \cite{Fornberg2015nodegen,Slak2019} and $N$ the total number of generated nodes. Larger $N$ suggests closer to optimal node placement.}
	\label{fig:comparisonunifnodes}
\end{figure}

As a first test case, nodes are generated in the unit square with constant exclusion radius $r = 0.025$. They are compared to nodes generated by the original Fornberg \& Flyer method \cite{Fornberg2015nodegen} and nodes generated by the recent Slak \& Kosec method \cite{Slak2019} as well as a Cartesian lattice. Note that both methods have a parameter $n$ to adjust, which corresponds to the number of sample points generated at each step. We use the recommended $n=5$ for Fornberg \& Flyer and $n=15$ for Slak \& Kosec. The three node sets can be seen in \cref{fig:comparisonunifnodes}. The optimal circle packing in the plane is hexagonal and visually one can see that the present method results in nodes most similar to this.
 
\begin{figure}[htbp]
	\centering
	\begin{subfigure}{0.45\textwidth}
		\centering
		\includegraphics[width = 0.85\textwidth]{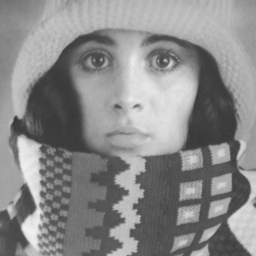}
		\caption{Original image}
		\label{fig:truioriginal}
	\end{subfigure}
	\begin{subfigure}{0.45\textwidth}
		\centering
		\includegraphics[width=0.875\textwidth]{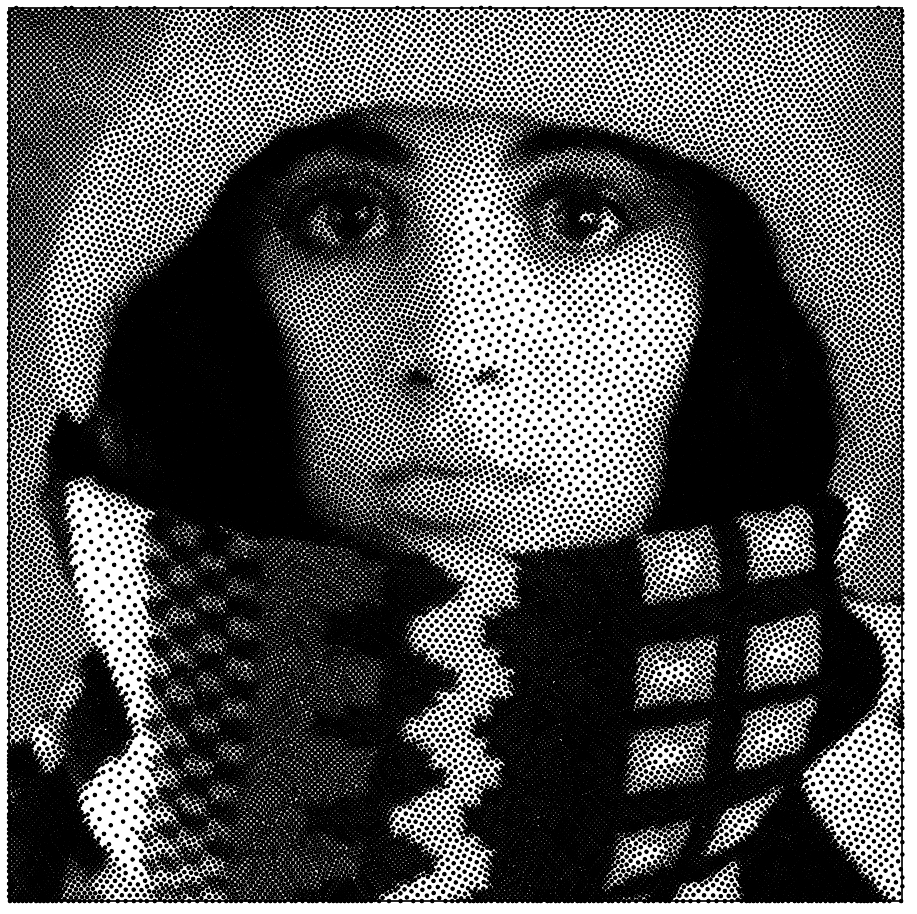}
		\caption{Dithered version}
		\label{fig:truidither}
	\end{subfigure}
	\caption{Test image `trui' and dithered version using the present method with a total of $N = 40,664$ nodes.}
	\label{fig:trui}
\end{figure}

\begin{figure}[htbp]
	\centering
	\begin{subfigure}{0.32\textwidth}
		\centering
		\includegraphics[width=0.9\textwidth]{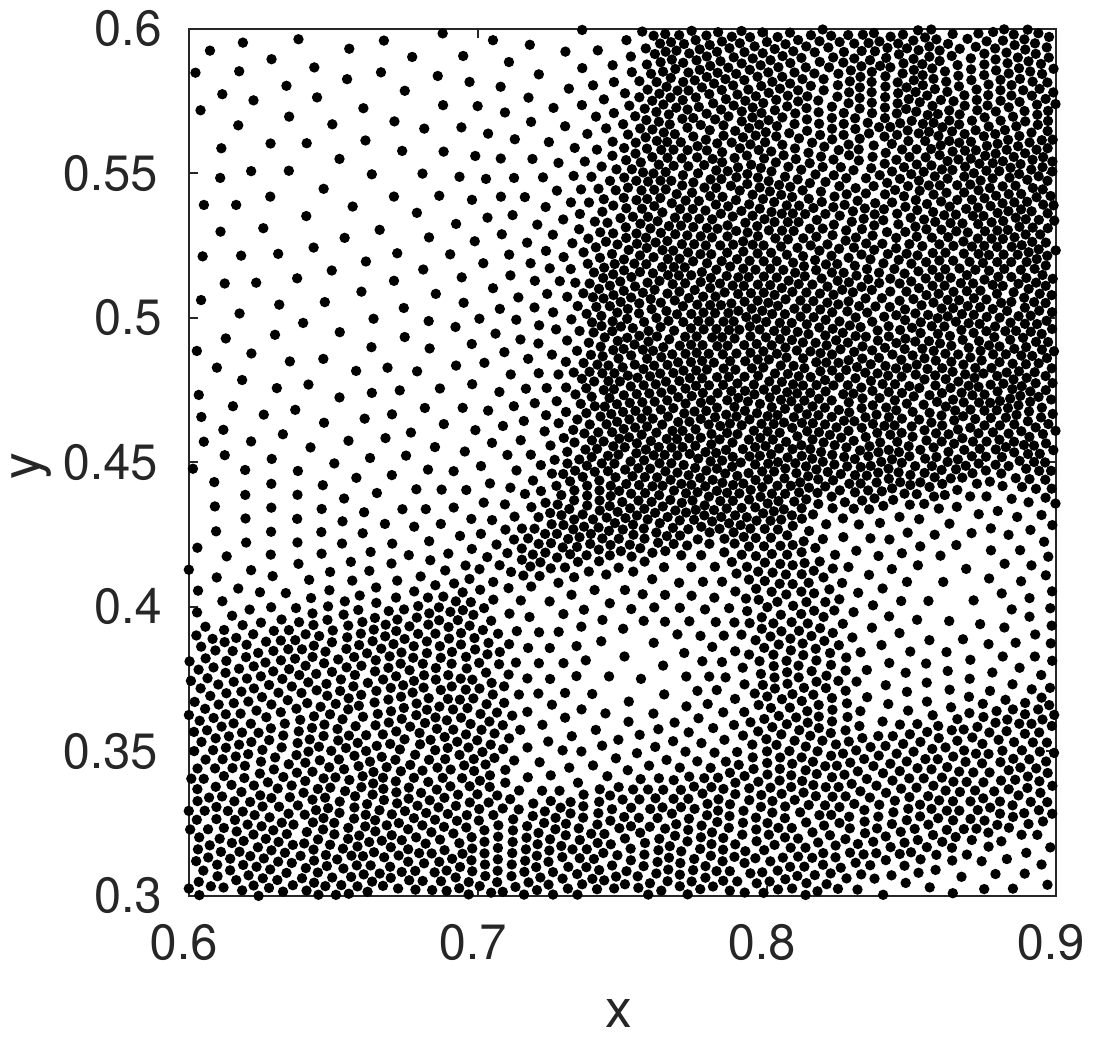}
		\caption{\centering 
			Fornberg \& Flyer, $N = 36,328$}
		\label{fig:FFvarnodes}
	\end{subfigure}
	\begin{subfigure}{0.32\textwidth}
		\centering
		\includegraphics[width=0.9\textwidth]{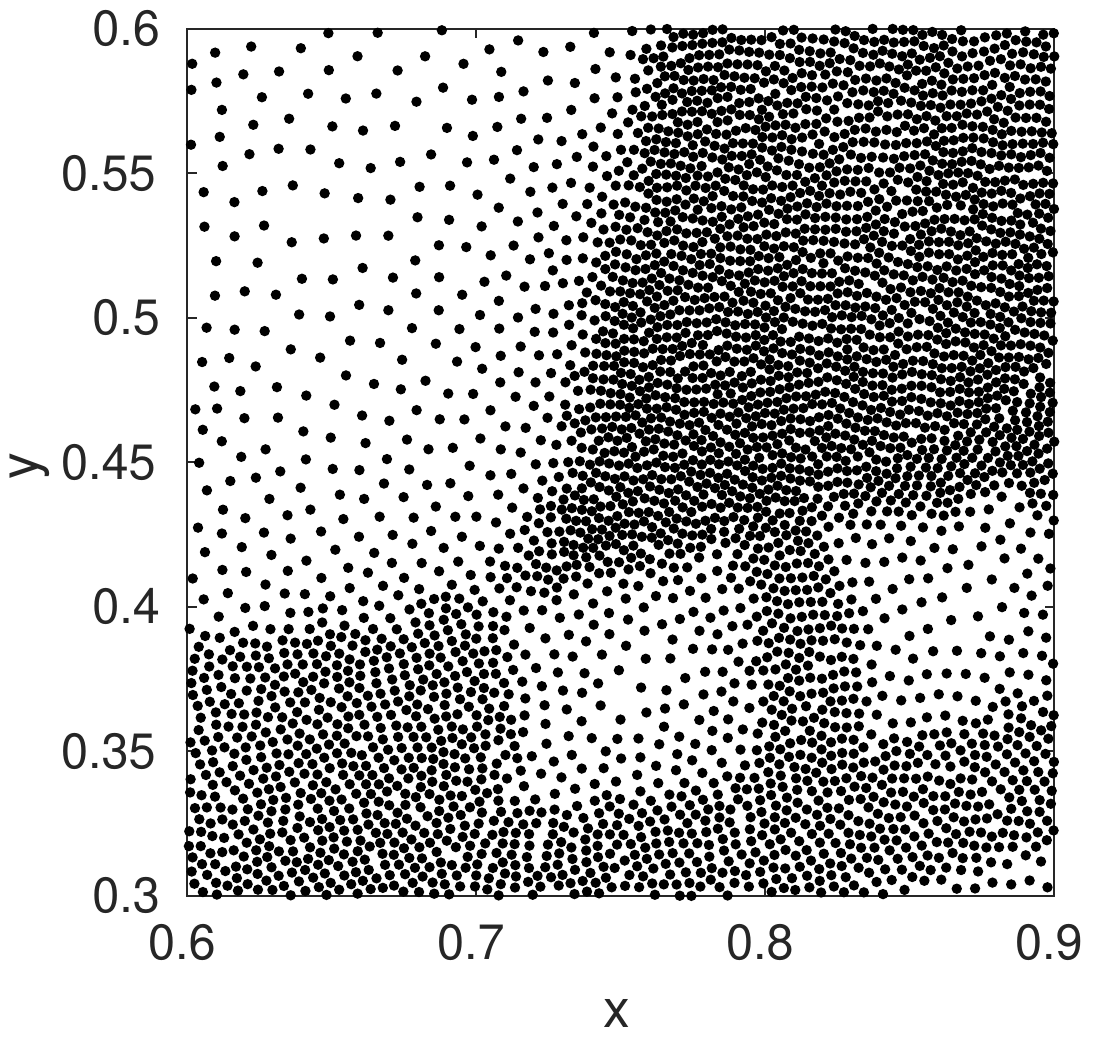}
		\caption{\centering 
			Slak \& Kosec, $N = 35,014$}
		\label{fig:SKvarnodes}
	\end{subfigure}
	\begin{subfigure}{0.32\textwidth}
		\centering
		\includegraphics[width=0.9\textwidth]{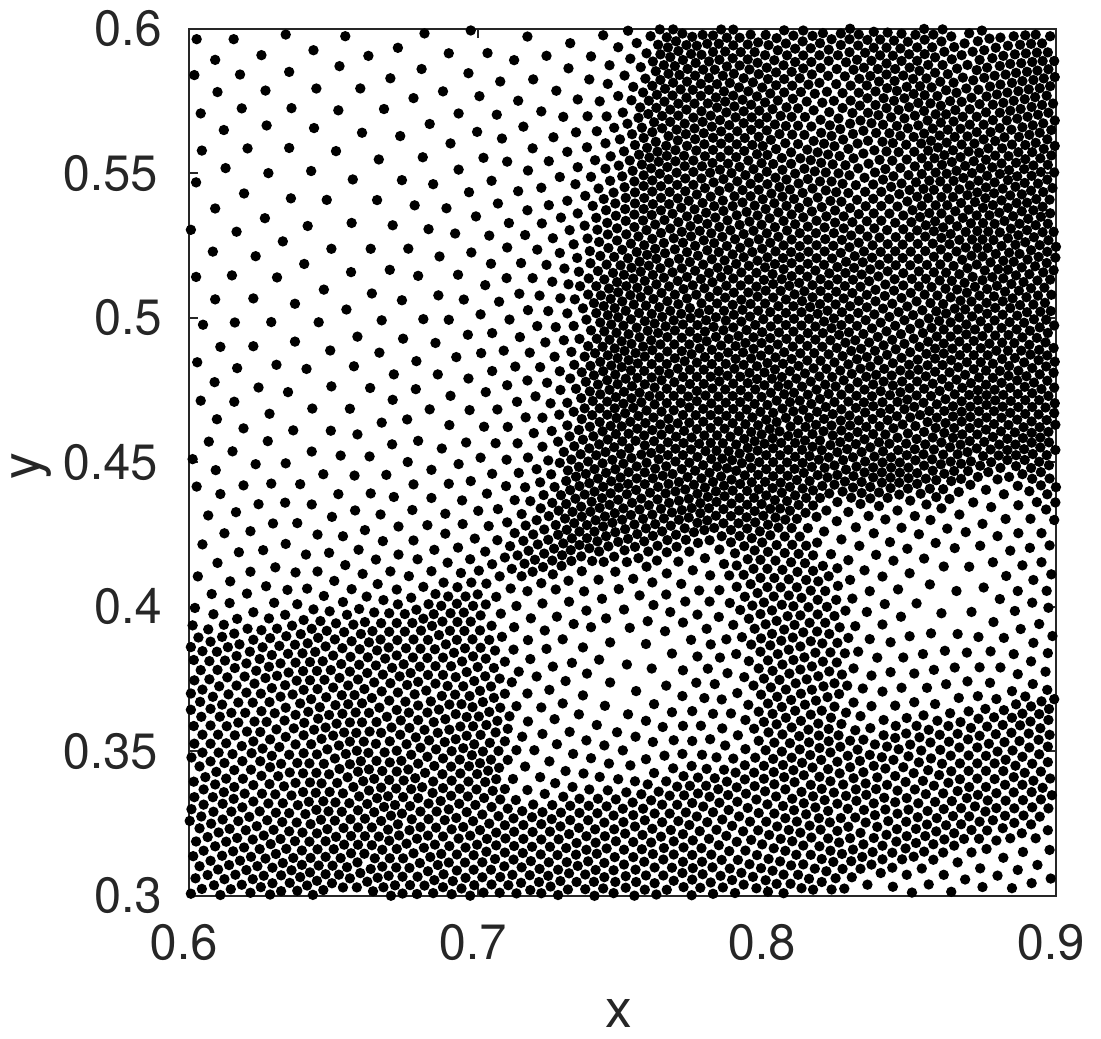}
		\caption{\centering
			 Present method, $N = 40,664$}
		\label{fig:newvarnodes}
	\end{subfigure}
	\caption{Enlargements of 2-D variable density node sets based off the `trui' image shown in \cref{fig:truioriginal}. $N$ is the number of nodes generated in the total dithered image.}
	\label{fig:comparisontruinodes}
\end{figure}

One of the advantages of the present method compared to a lattice structure or Halton node set is the ability to handle prescribed variable density functions. To demonstrate the ability to generate locally quasi-uniform nodes of highly variable density with sharp gradients, we consider the common test case for image rendering found online as `trui.png'. The radial exclusion function is based off the gray-scale information of the image so that more nodes are placed in darker areas. \cref{fig:trui} depicts the original image and the resulting dithered nodes using the present method, while \cref{fig:comparisontruinodes} shows a close-up comparing the three algorithms. Note that the same radial exclusion function results in different numbers of generated nodes for each method. The present algorithm is able to achieve the highest density while still satisfying the minimum spacing constraint of the function. Qualitatively, the nodes also look the most locally regular. 

\cref{table:2d_fixedspacing} compares the quality metrics defined in \cref{sec:nodesetquality} for the uniform node sets and \cref{table:2d_variabledensity} compares the dithered node sets. Both the mesh ratio and packing density are unit free metrics that don't depend on the number of nodes being placed. The present method has the highest number of nodes and smallest mesh ratio $\gamma$. The packing density gives a way to compare to the optimal hexagonal circle packing density $\pi\sqrt{3}/6 \approx 0.9069$. The present method is closest to this optimal density. Note that for the uniform case, $\rho$ should be close to half of the prescribed spacing $r=0.025$. For the variable density case, $\rho$ is not included since as the exclusion radius is varying the ratio $\gamma$ is more descriptive. 

\begin{table}[htbp]
	\centering
	\begin{tabular}{lcccc}
		\hline
		\noalign{\smallskip}
		\multicolumn{5}{c}{Uniform Density} \\ 
		\hline
		\noalign{\smallskip}
		& & Packing &  Covering Radius & Mesh Ratio \\ 
		Method & N & Density & $\rho$ & $\gamma$ \\ 
		\hline
		\noalign{\smallskip}
		Fornberg \& Flyer & 1492 & 0.72 & 0.0216 & 0.746 \\
		Slak \& Kosec & 1429 & 0.70 & 0.0253 & 1.010 \\
		Cartesian & 1681 & 0.79 & 0.0177 & 0.707 \\
		Present method & 1698 & 0.83 & 0.0180 & 0.685 \\ 
		\hline
		\noalign{\smallskip}
		Theoretical limit & - & 0.91 & 0.0125 & 0.5
	\end{tabular}	
	\caption{A comparison of node quality metrics on 2-D uniform node sets with spacing $r=0.025$. Larger $N$, larger packing density, smaller $\rho$ and smaller $\gamma$ are better.}
	\label{table:2d_fixedspacing}
\end{table}

\begin{table}[htbp]
	\centering
	\begin{tabular}{lcc}
		\hline
		\noalign{\smallskip}
		\multicolumn{3}{c}{Variable Density} \\ 
		\hline
		\noalign{\smallskip}
		& & Mesh Ratio \\ 
		Method & N & $\gamma$ \\ 
		\hline
		\noalign{\smallskip}
		Fornberg \& Flyer & 36,328 & 0.756 \\
		Slak \& Kosec & 35,014 & 0.780 \\
		Present method & 40,663 & 0.637 \\ 
		\hline
	\end{tabular}	
	\caption{A comparison of node quality metrics on variable density node sets generated from the `trui' test case. Larger $N$, larger packing density, smaller $\rho$ and smaller $\gamma$ are better.}
	\label{table:2d_variabledensity}
\end{table}

Looking at the mesh norm $\gamma$ gives an idea of whether the nodes are well-spaced and quasi-uniform. However, it is not a perfect metric. 
These metrics only take into account the maximum gap and minimum distance between neighbors, not the distribution of gaps over the whole node set. 
A Cartesian grid will have the same gap and minimum distance over the whole node set while nodes generated by the present method have smaller gaps overall, which is why $N$ is closer to the maximal packing, but may have a few nodes with larger gaps than the grid. This is why the Cartesian grid has a smaller covering radius than the present method. In fact Cartesian grids are non-optimal for RBF-FD due to poor conditioning and accuracy issues \cite{FF2015primer}, which are investigated in \cref{sec:RBFmethods}. 

To get a better idea of the variation over the node set, local regularity can be observed from the distribution of distance to the nearest $k$ neighbors $\delta_{i,j}$, $i=1,2,...k$ for each node $p_j$.  \cref{table:2dneighborcomparison} compares statistics based on 6 nearest neighbors for the uniform node sets. We use $k=6$ based on hexagonal circle packing. Here, the present method gives the closest $\bar{\delta_j}$ to the prescribed $r=0.025$ with a small standard deviation and mean range($\delta_{i,j}$). 

\begin{table}[htbp]
	\centering
	\begin{tabular}{lccc}
		\hline
		\noalign{\smallskip}
		Method & mean $\bar{\delta_j}$ &  std $\bar{\delta_j}$ & mean range($\delta_{i,j}$) \\ 
		\hline
		\noalign{\smallskip}
		Fornberg \& Flyer & 0.0291 & 0.0016 & 0.0114 \\
		Slak \& Kosec & 0.0299 & 0.0018 & 0.0127 \\
		Cartesian & 0.0285 & 3.3e-16 & 0.0104 \\
		Present method & 0.0270 & 0.0007 & 0.0080 
	\end{tabular}	
	\caption{A comparison of distance to nearest 6 neighbors for uniform 2-D nodes. Mean $\bar{\delta_j}$ should be close to prescribed $r=0.025$ while std $\bar{\delta_j}$ and mean range($\delta_{i,j}$) should be small.}
	\label{table:2dneighborcomparison}
\end{table}

\begin{figure}[htbp]
	\centering
	\begin{subfigure}[]{0.34\textwidth}
		\centering
		\includegraphics[width=0.9\textwidth]{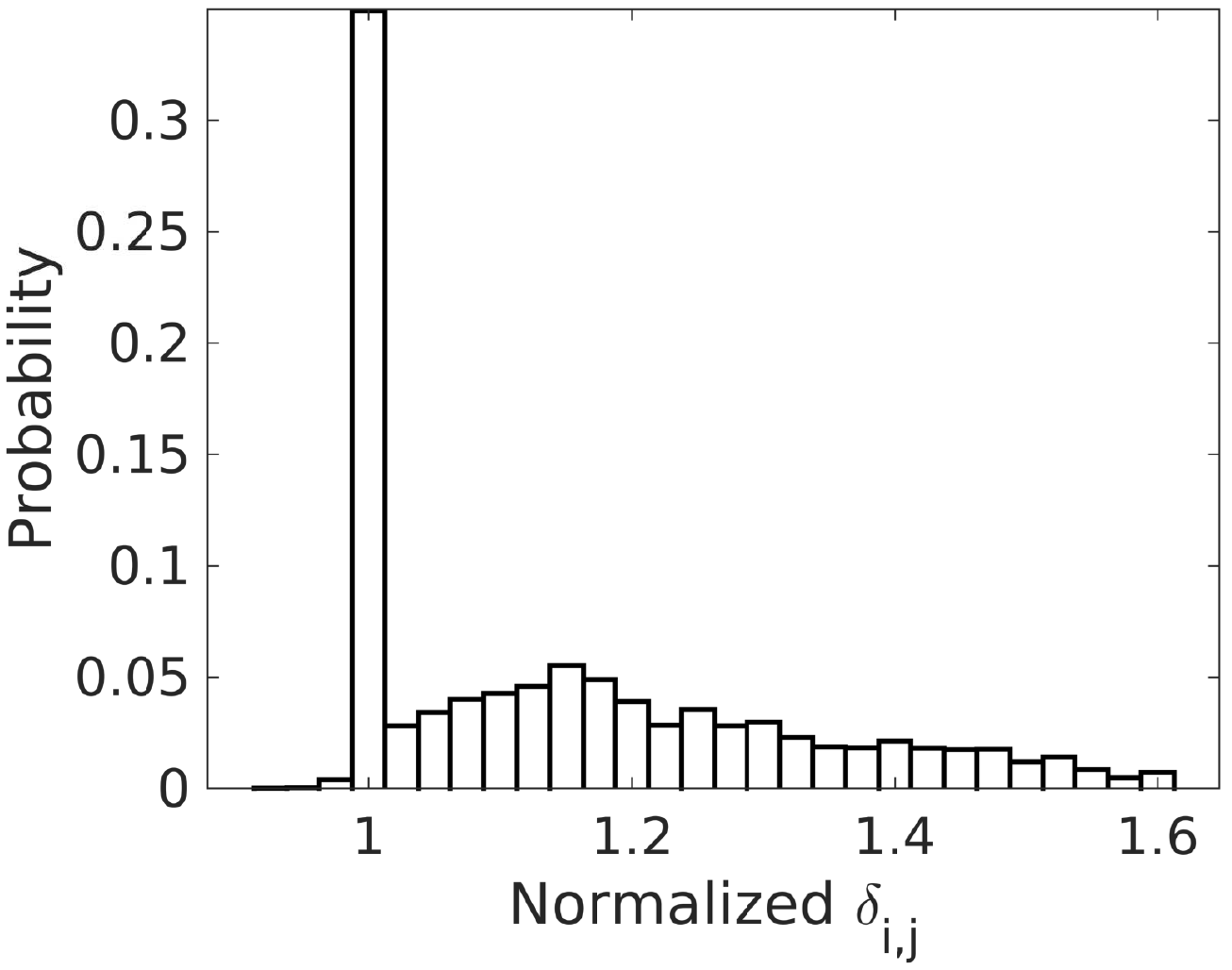}
		\caption{Fornberg \& Flyer}
	\end{subfigure}
	\begin{subfigure}[]{0.34\textwidth}
		\centering
		\includegraphics[width=0.9\textwidth]{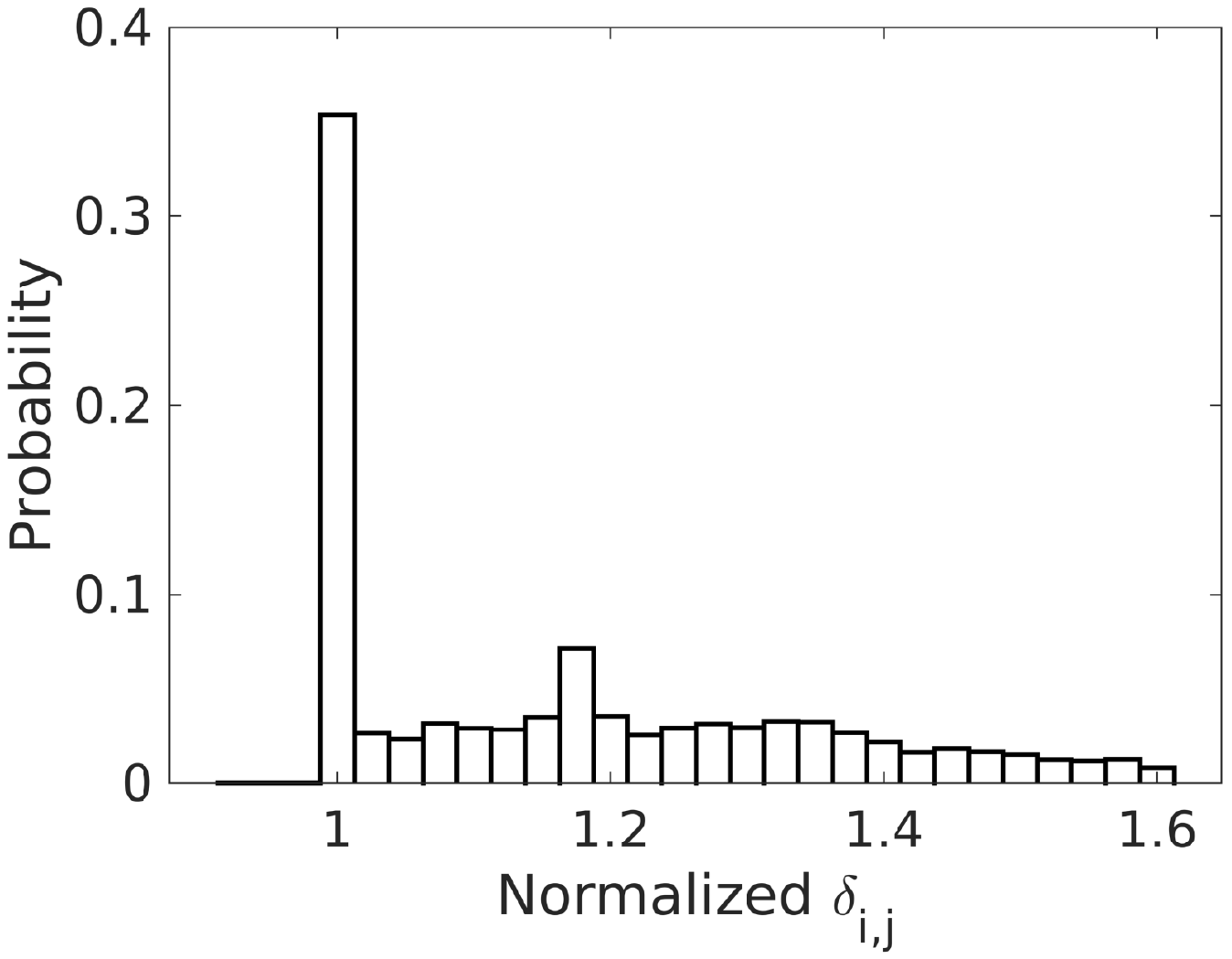}
		\caption{Slak \& Kosec}
	\end{subfigure}
	\begin{subfigure}[]{0.34\textwidth}
		\centering
		\includegraphics[width=0.9\textwidth]{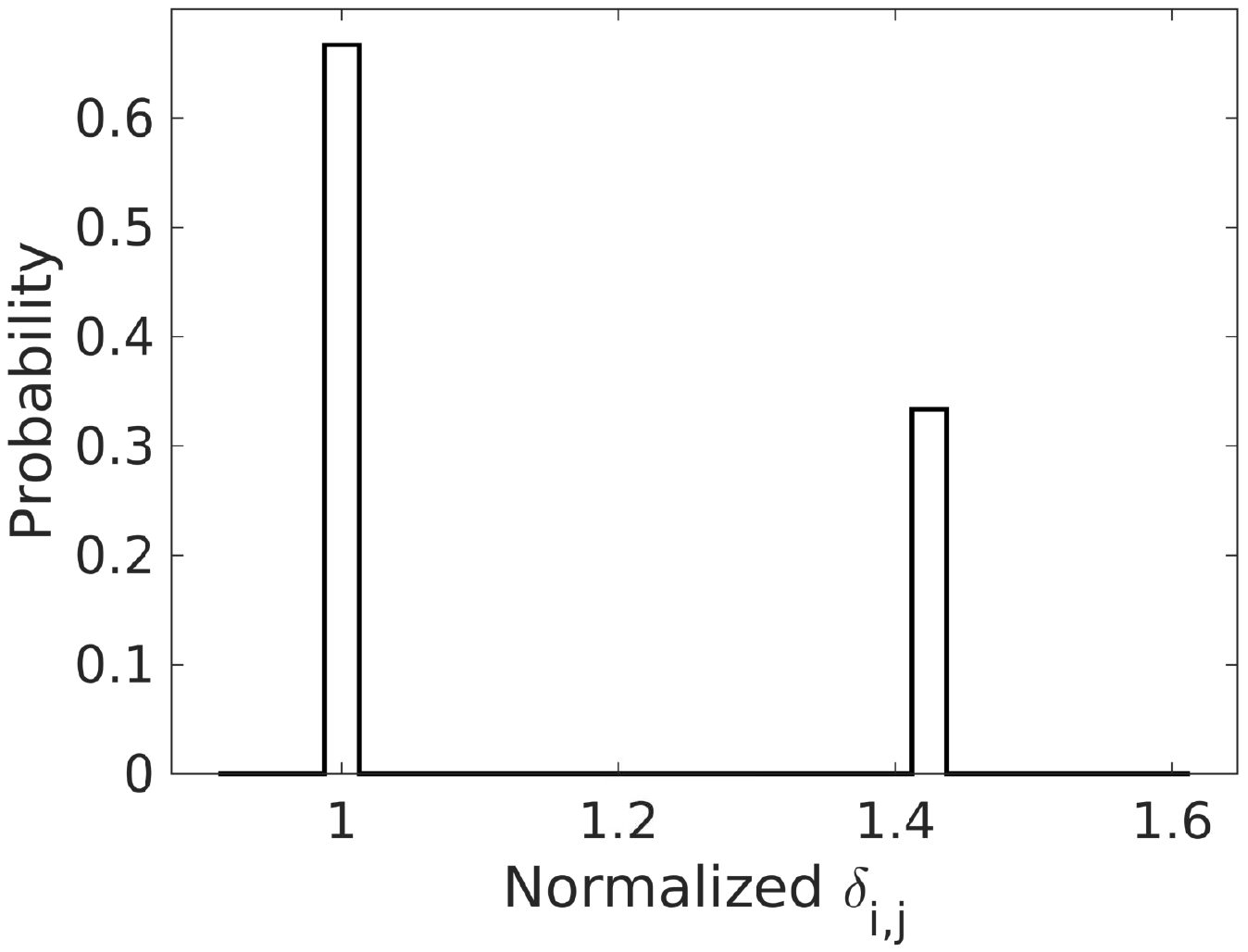}
		\caption{Cartesian lattice}
	\end{subfigure}
	\begin{subfigure}[]{0.34\textwidth}
		\centering
		\includegraphics[width=0.9\textwidth]{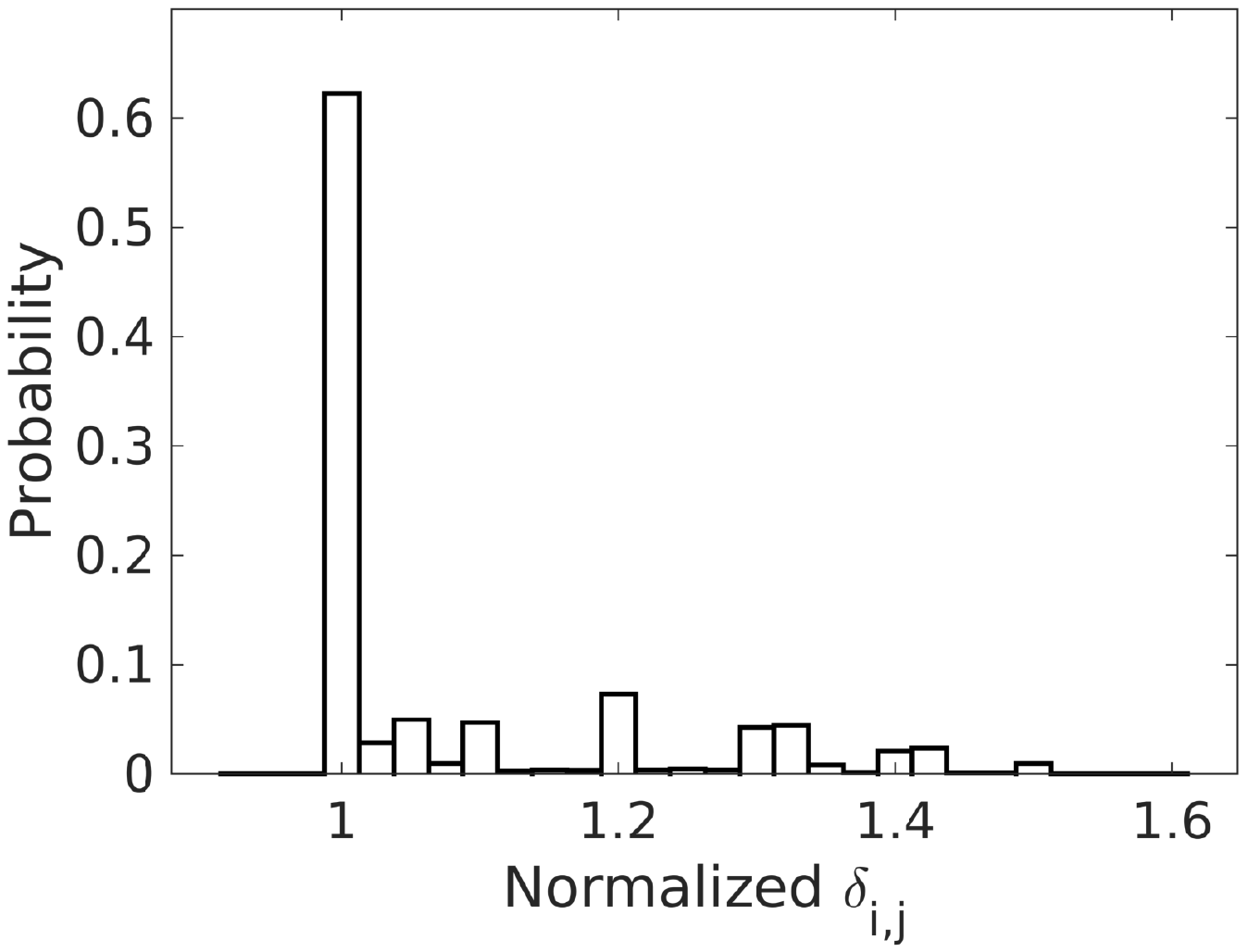}
		\caption{Present method}
	\end{subfigure}
	\caption{Distribution of distance to 6 nearest neighbors for the uniform 2-D node sets. Distance to neighbors $\delta_{i,j}$ is normalized by the exclusion radius $r = 0.025$
	and the number of counts in each bin is normalized by the total number of counts
		.}
	\label{fig:2duniformnodehist}
\end{figure}

One way to visualize this distribution of nearest neighbors is through a histogram plot, as seen in \cite{,Shankar2018,Vlasiuk2018}. The distance to nearest neighbor can be scaled by the exclusion radius function so a sharp peak is expected around 1 with some spread to the right. This can be seen in the histograms in \cref{fig:2duniformnodehist}. The 6 neighbors in the Cartesian lattice are fixed at one of two distances as expected in a grid lattice. More neighbors are at the prescribed distance in the present method than Slak \& Kosec.  

\subsection{Nodes in 3-D}
\label{sec:3Dresults}
Moving on to the 3-D case, nodes were generated in the unit cube. A uniform node set with prescribed exclusion radius $r=0.05$ is compared to a Cartesian lattice and a node set generated with the method from \cite{Slak2019} in \cref{table:3dmetriccomparison}. Note there is no comparison to \cite{Fornberg2015nodegen} as this method does not generalize to higher dimensions. As in the 2-D case, $\rho$ is smaller for a Cartesian grid than for the nodes generated by the present method. On all other metrics the present method performs best. 

Optimal packing density in higher than 2 dimensions is a classic challenging mathematical problem \cite{Conway1993}. The optimal 3-D packing density is achieved by a family of close packed lattices, which have a packing density of $\pi/3\sqrt{2} \approx 0.74$. For the present method, background density was increased to a factor of 100 to determine the average packing density. Using the usual factor of 10 would decrease cost and only compromise on quality by about $5\%$.

\begin{table}[htbp]
	\centering
	\begin{tabular}{lcccc}
		\hline
		\noalign{\smallskip}
		Method & N & Packing Density &Covering Radius $\rho$ & Mesh Ratio $\gamma$ \\ 
		\hline
		\noalign{\smallskip}
		Slak \& Kosec & 7128 & 0.46 & 0.0537 & 1.073 \\
		Cartesian & 9261 & 0.52 & 0.0433 & 0.866 \\
		Present method & 9998 & 0.61 & 0.0447 & 0.887 \\ 
		\hline
		\noalign{\smallskip}
		Theoretical limit & - & 0.74 & 0.0250 & 0.500
	\end{tabular}	
	\caption{A comparison of node quality metrics on 3-D uniform node sets with spacing $r=0.05$. Larger $N$, larger packing density, smaller $\rho$ and smaller $\gamma$ are better.}
	\label{table:3dmetriccomparison}
\end{table}

\cref{table:3dneighborcomparison} compares statistics based on 12 nearest neighbors for the same 3 uniform node sets. We use $k=12$ based on dense sphere packings like cubic close packing and hexagonal close packing where each sphere is surrounded by 12 others. Here, the present method gives the closest $\bar{\delta_j}$ to the prescribed $r=0.05$ with a small standard deviation and mean range($\delta_{i,j}$). 
	
\begin{table}[htbp]
	\centering
	\begin{tabular}{lccc}
		\hline
		\noalign{\smallskip}
		Method & mean $\bar{\delta_j}$ &  std $\bar{\delta_j}$ & mean range($\delta_{i,j}$) \\ 
		\hline
		\noalign{\smallskip}
		Slak \& Kosec & 0.0608 & 0.0017 & 0.0254 \\
		Cartesian & 0.0604 & 3.6e-15 & 0.0207 \\
		Present method & 0.0546 & 8.8e-4 & 0.0163 
	\end{tabular}	
	\caption{A comparison of distance to nearest 12 neighbors for uniform 3-D nodes. Mean $\bar{\delta_j}$ should be close to prescribed $r=0.05$ while std $\bar{\delta_j}$ and mean range($\delta_{i,j}$) should be small.}
	\label{table:3dneighborcomparison}
\end{table}
	
\begin{figure}[!htbp]
	\centering
	\begin{subfigure}[]{0.32\textwidth}
		\centering
		\includegraphics[width=0.9\textwidth]{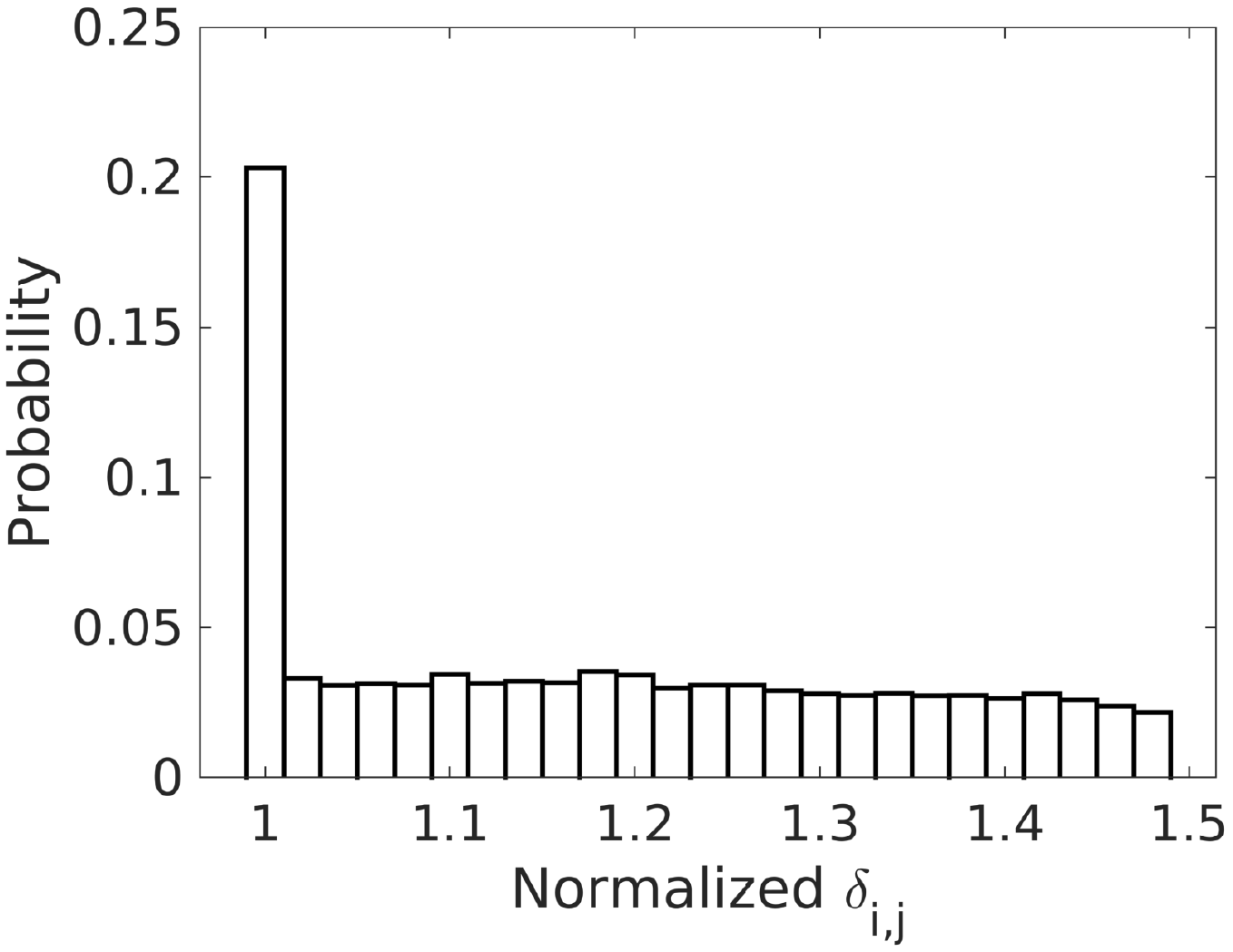}
		\caption{Slak \& Kosec}
	\end{subfigure}
	\begin{subfigure}[]{0.32\textwidth}
		\centering
		\includegraphics[width=0.9\textwidth]{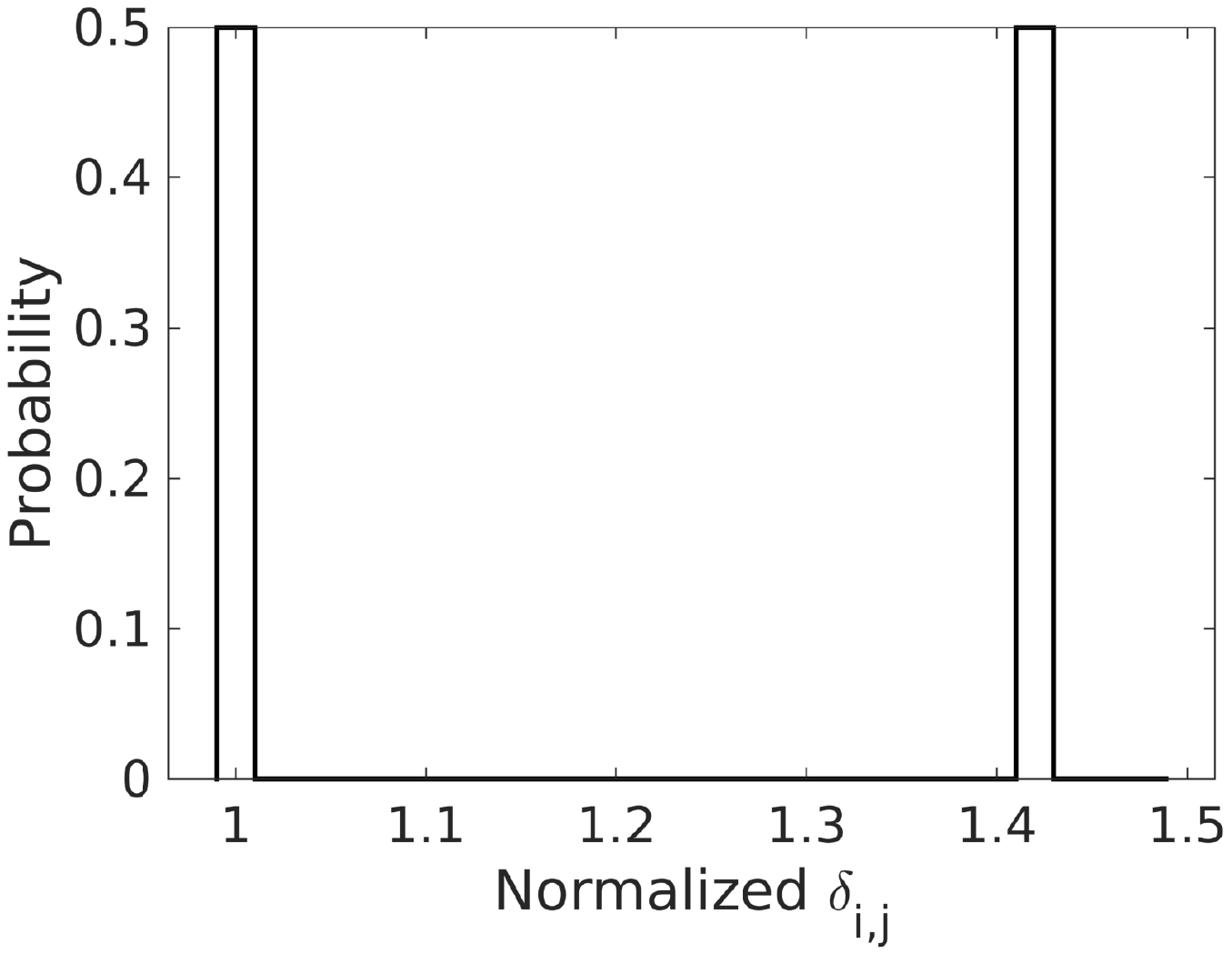}
		\caption{Cartesian lattice}
	\end{subfigure}
	\begin{subfigure}[]{0.32\textwidth}
		\centering
		\includegraphics[width=0.9\textwidth]{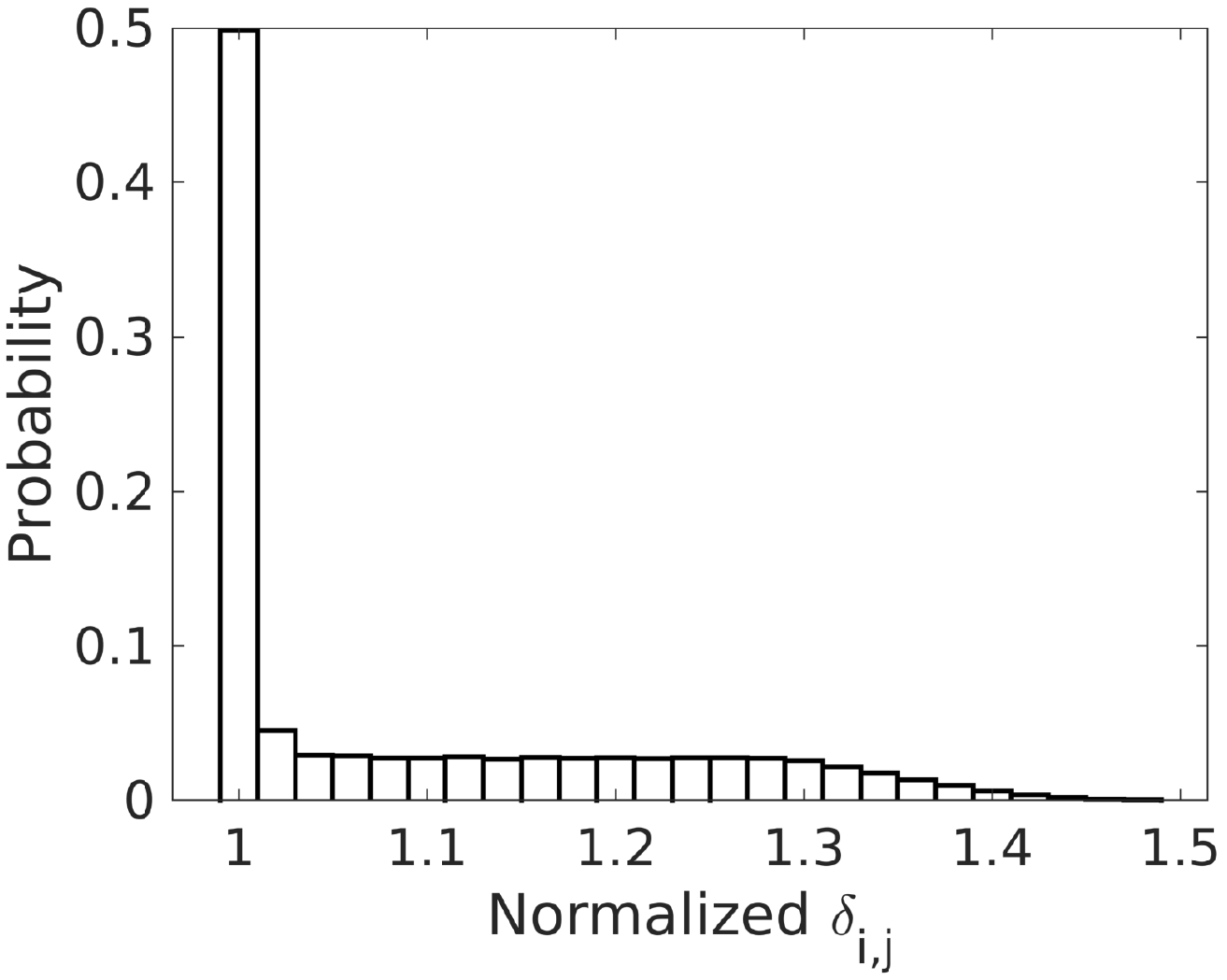}
		\caption{Present method}
	\end{subfigure}
	\caption{Distribution of distance to 12 nearest neighbors for the uniform 3-D node sets. Distance to neighbors $\delta_{i,j}$ is normalized by the exclusion radius $r = 0.05$
	.}
	\label{fig:uniformnodehist}
\end{figure}

Histogram plots in \cref{fig:uniformnodehist} show the distribution of distance to the nearest 12 neighbors. Again the neighbors in the Cartesian lattice are fixed at one of two distances as expected in a grid lattice. More neighbors are at the prescribed distance in the present method than Slak \& Kosec, with a smoother tail.  

\subsection{Execution Time}
\begin{figure}[htbp]
	\centering
	\begin{subfigure}[]{0.45\textwidth}
		\centering
		\includegraphics[width=0.8\textwidth]{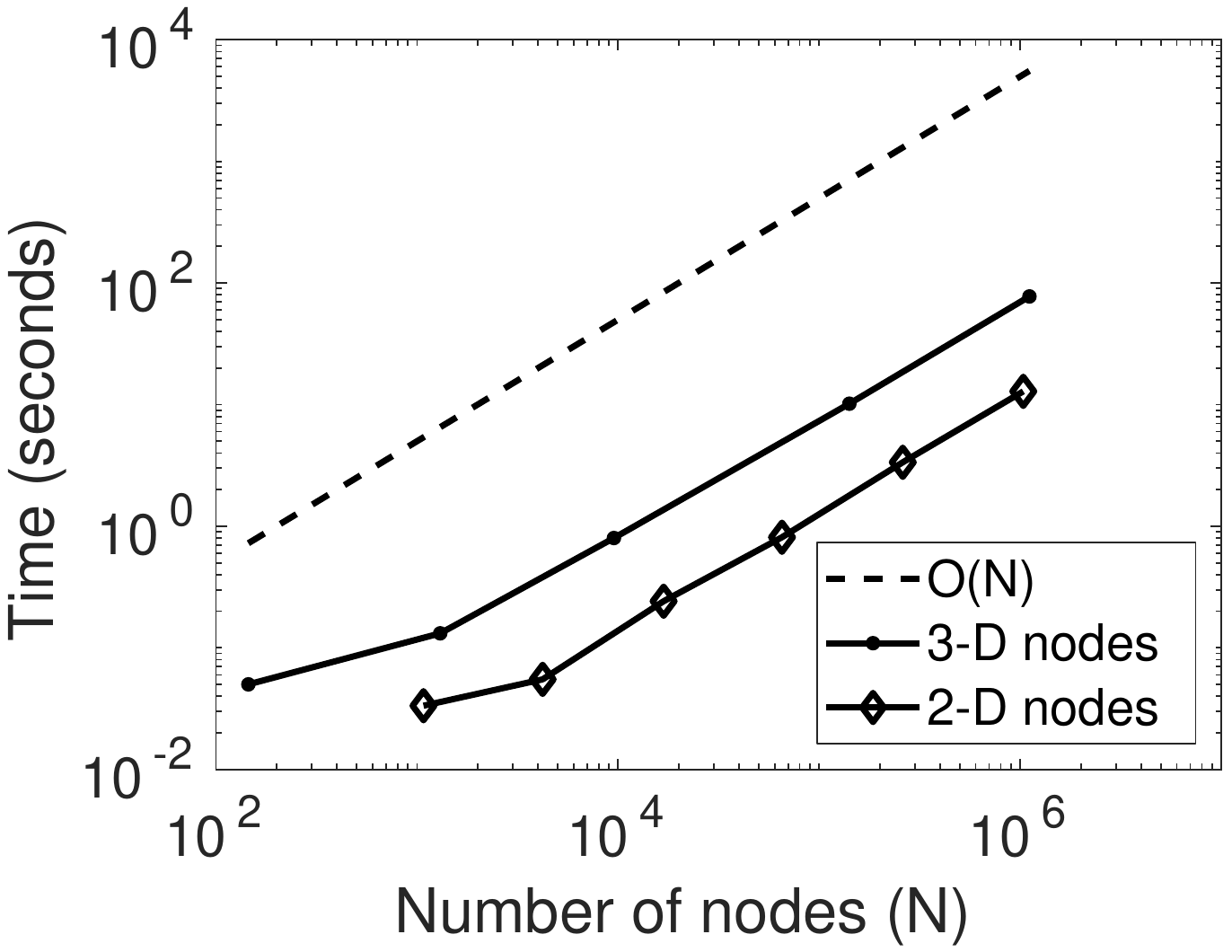}
		\caption{Uniform density}
		\label{fig:uniformnodesconv}
	\end{subfigure}
	\begin{subfigure}[]{0.45\textwidth}
		\centering
		\includegraphics[width=0.8\textwidth]{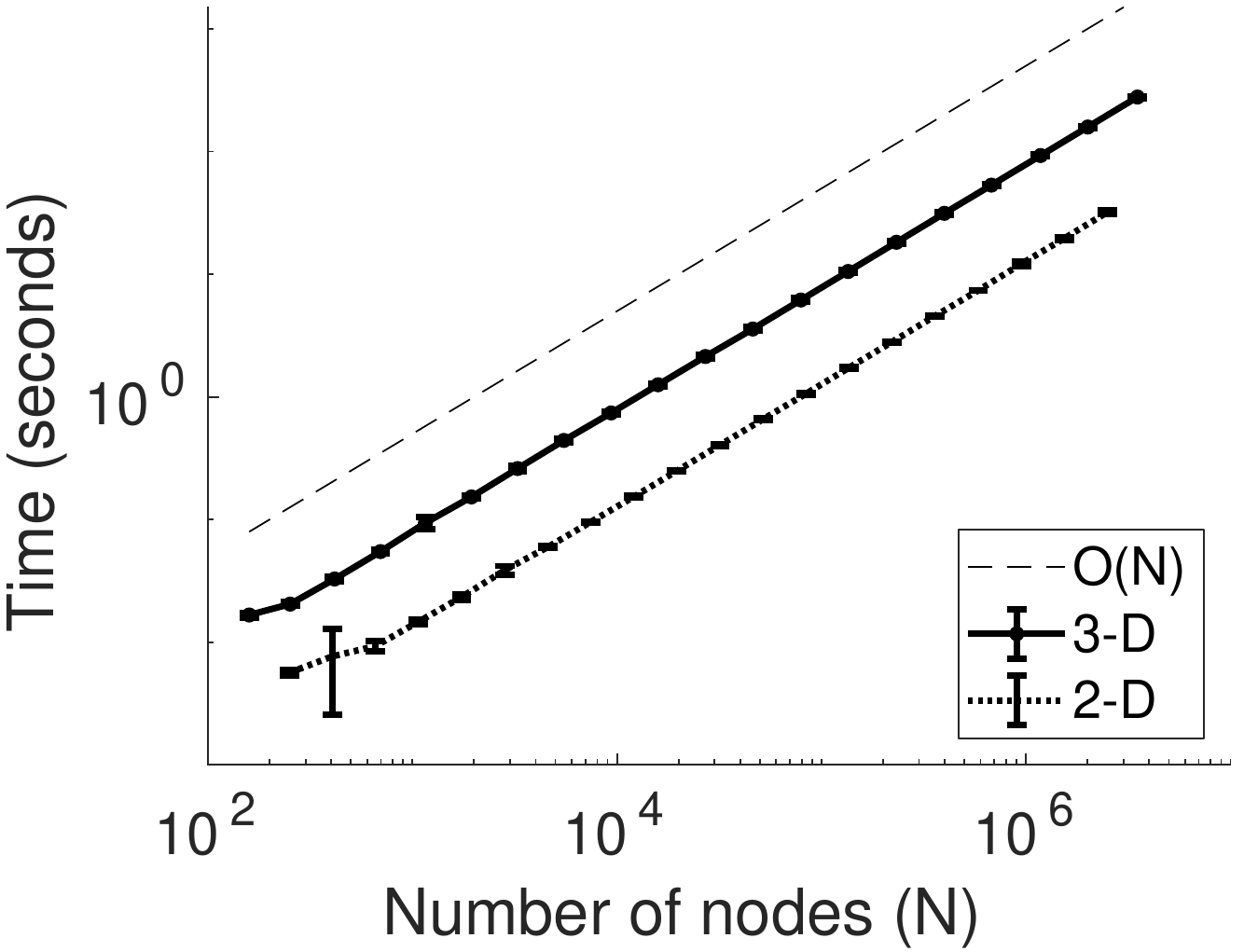}
		\caption{Variable density}
	\end{subfigure}
	\caption{Cost of node generation in the 2-D unit square and 3-D unit cube with no boundary using a MATLAB implementation on a 6-core \textit{Intel i7-8750H} CPU. An average is taken over 10 tests per each number of nodes, and error bars denote standard deviation from the mean.}
	\label{fig:nodegenerationcost}
	\end{figure}

We investigate the time complexity of the present method through numerical experiments. Node generation cost is expected to scale with number of nodes $N$ placed for a fixed background grid density factor. In \cref{fig:nodegenerationcost} we observe $O(N)$ convergence for both uniform density and variable density node sets in 2-D and 3-D. 

There are two \texttt{while} loops in the given algorithm. The outer loop runs until the nodes are out of the bounding box. Since the algorithm satisfies minimum spacing between nodes, an upper bound on the number of nodes to be placed in bounding box $B$ of dimension $d$, and thus the number of loop iterations, is given by
\begin{equation}
	N_{max} < \frac{\text{Vol}(B)}{\frac{\pi^{d/2}}{\Gamma\left(\frac{d}{2}+1\right)}\left(\frac{\text{min}(r)}{2}\right)^d}. \label{eqn:Nbound}
\end{equation}
For a bounded box volume and $\text{min}(r)\geq \delta > 0$, this is a finite bound. 

The inner loop is an iteration to find a local minimum. In the worst case scenario, this minimum search will continue until the global minimum is found. The global minimum exists for a bounded box, since the projection onto the plane orthogonal to the $z$ direction is also a finite area which we denote $bb$. Thus the maximum number of iterations for node $p$ is bounded as
\begin{equation}
	N_{iter} \leq \frac{bb -(2r(p))^{d-1}}{r(p)^{d-1}},
\end{equation}
which is also finite for an exclusion radius function $r(p)\geq \delta > 0$. 

In practice, the number of iterations to find a local minimum is significantly less than this upper bound. For the 2-D uniform node set shown in \cref{fig:newnodes} the average number of iterations per node is $1.72$, while for the 2-D \texttt{trui} node set the average is $1.67$. For the 3-D uniform node set the average number of iterations is $2.08$. Although the maximum number of iterations does increase as the density increases, the average remains around 2 in all the experiments detailed in this work. 

\subsection{Direction dependence correction results}
\begin{figure}[htbp]
	\centering
	\includegraphics[width=0.35\textwidth]{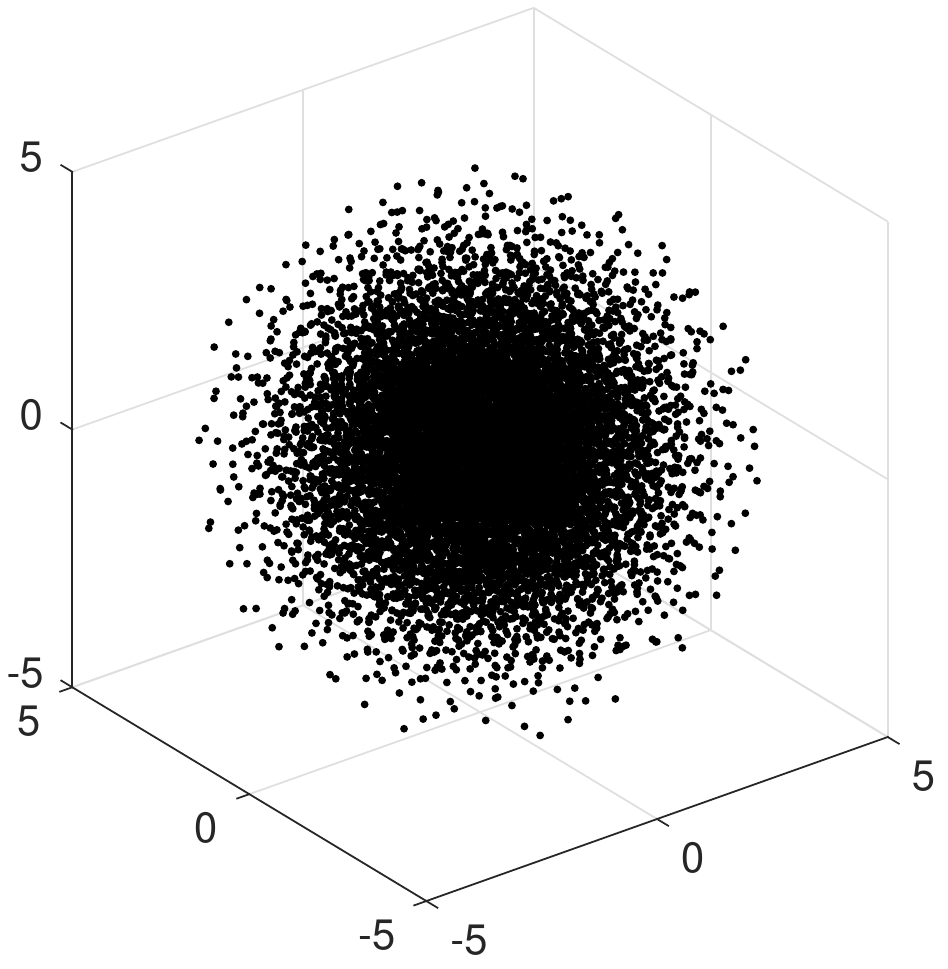}
	\caption{Radially varying node set generated using exclusion radius function \eqref{eqn:radialfunction3}.}
	\label{fig:radialnodes}
\end{figure}

The original algorithm is compared to the correction method described in \cref{sec:directioncorrection} and the radially built method described in \cref{sec:Spherical} for a radially varying node density. As we have already compared the present algorithm to other recent methods in \cref{sec:2Dresults} and \cref{sec:3Dresults}, we only compare to the present algorithm in this section. The exclusion radius function
\begin{equation}
	r(R) = Ce^{\epsilon R^2} \label{eqn:radialfunction3}
\end{equation}
is used as a test case, where $R = \sqrt{x^2 + y^2 + z^2}$ is the distance to the origin, and $C$ and $\epsilon$ are parameters that change the shape of the function. \cref{fig:radialnodes} shows a node set using $C = 4/21$ and $\epsilon = 1/15$. 

Two radially built variations are investigated. Variation (1) uses an exclusion radius based on previously placed nodes, while variation (2) uses the maximum exclusion radius between the current node and the previously placed nodes. 

\begin{figure}[htbp]
	\centering
	\begin{subfigure}{0.32\textwidth}
		\centering
		\includegraphics[width=0.9\textwidth]{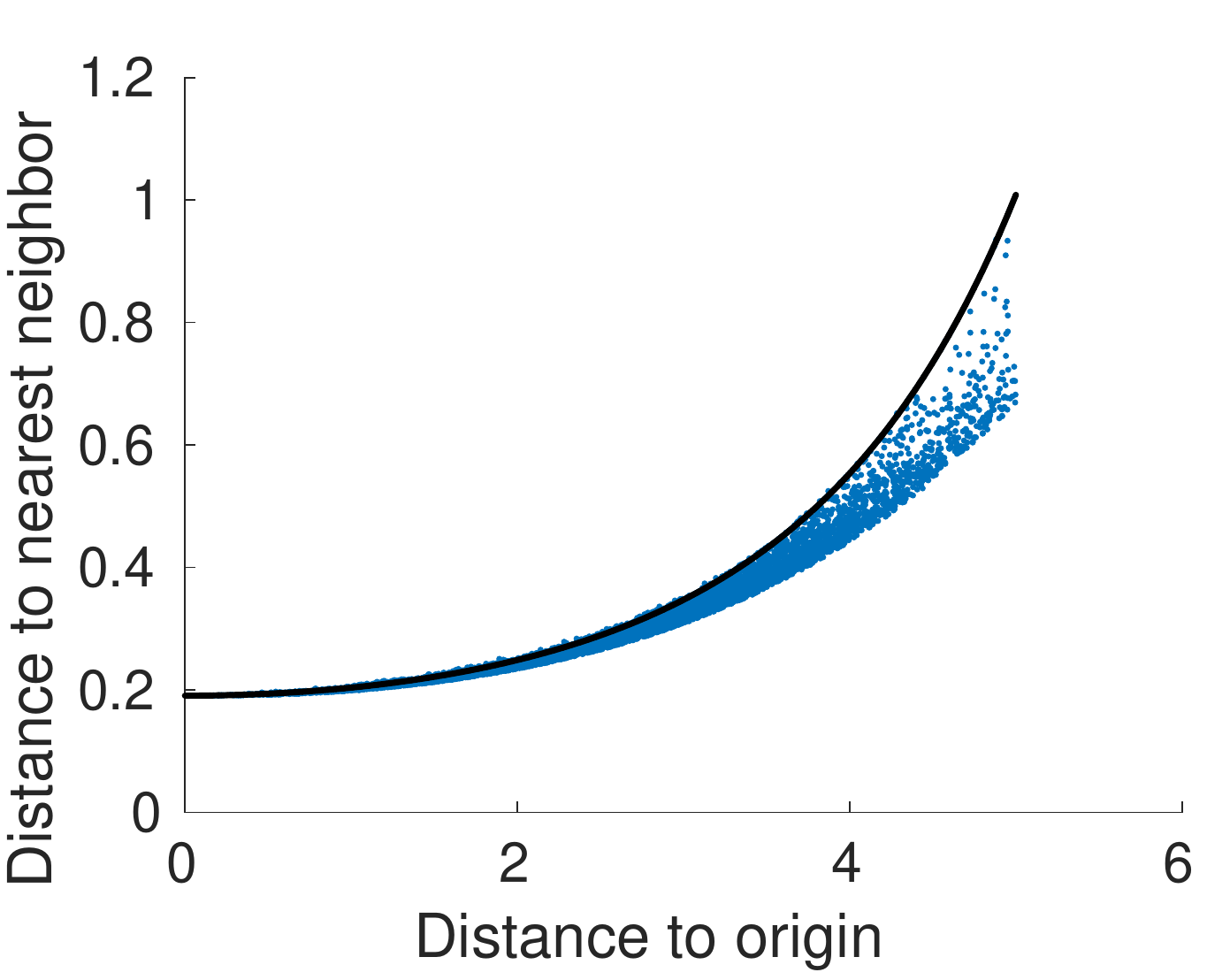}
		\caption{Original}
	\end{subfigure}
	\begin{subfigure}{0.32\textwidth}
		\centering
		\includegraphics[width=0.9\textwidth]{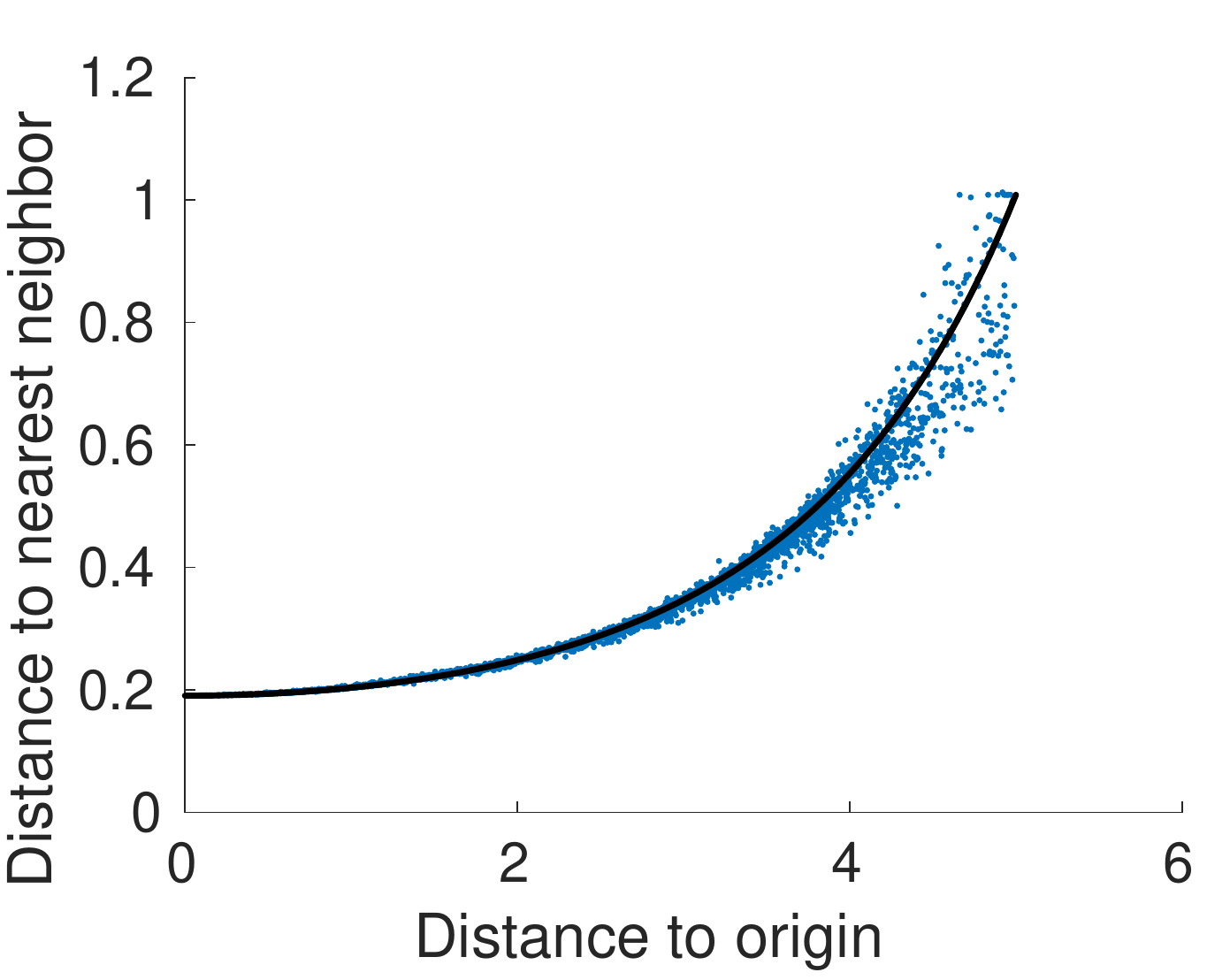}
		\caption{Corrected}
	\end{subfigure}
	\begin{subfigure}{0.32\textwidth}
		\centering
		\includegraphics[width=0.9\textwidth]{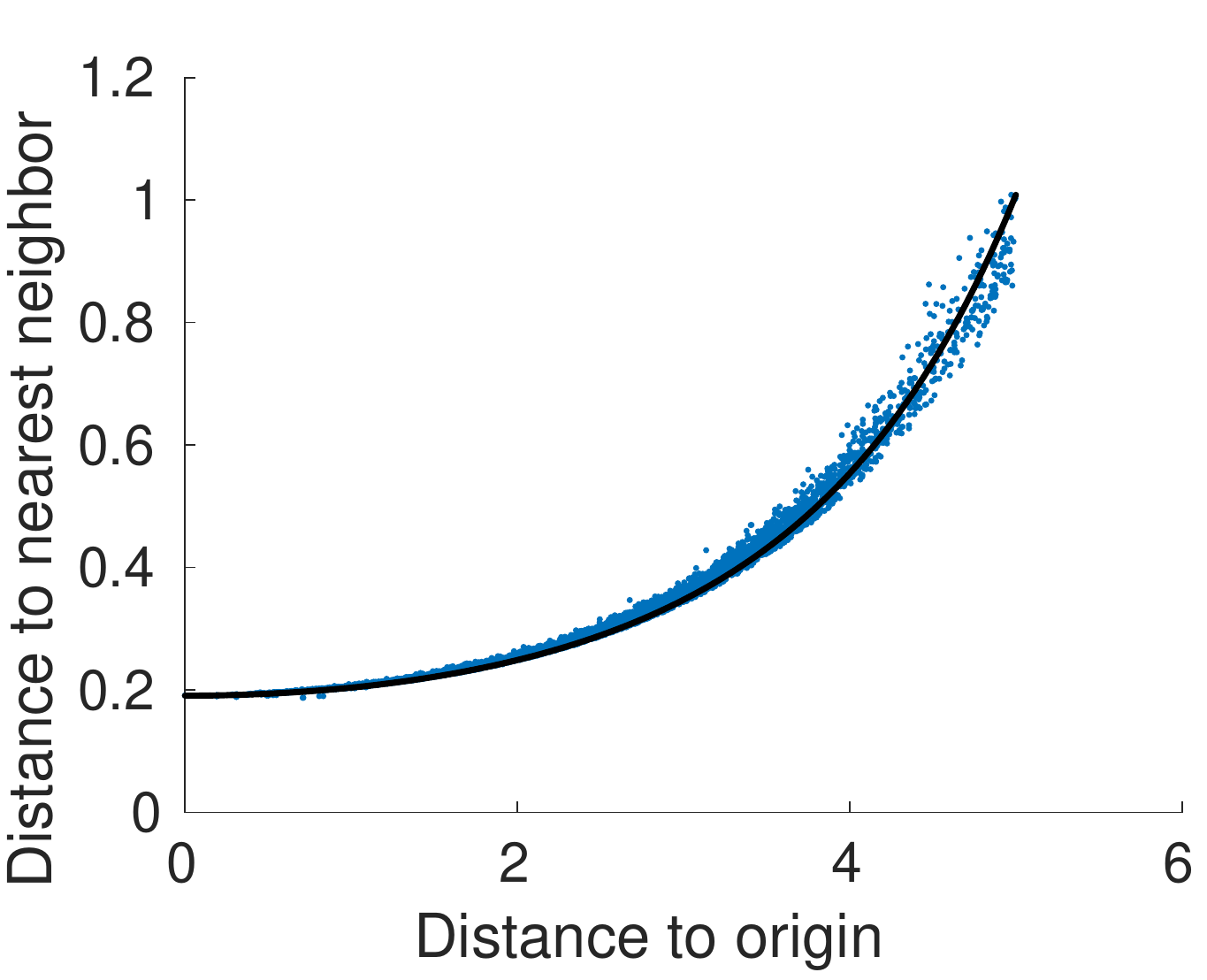}
		\caption{Radially built (1)}
	\end{subfigure}
	\caption{Scatter plots of the distance to nearest neighbor as a function of distance from the origin. The line through the data is the prescribed exclusion radius.}
	\label{fig:comparisonindistfromorig}
\end{figure}

First, the distance to nearest neighbor is compared to the prescribed exclusion radius function, both as a function of distance to the origin, in \cref{fig:comparisonindistfromorig}. The corrected version more closely aligns to the exclusion radius function than the original algorithm, but the radially built node set does even better. 

\begin{figure}[htbp]
	\centering
	\begin{subfigure}{0.325\textwidth}
		\centering
		\includegraphics[width=0.95\textwidth]{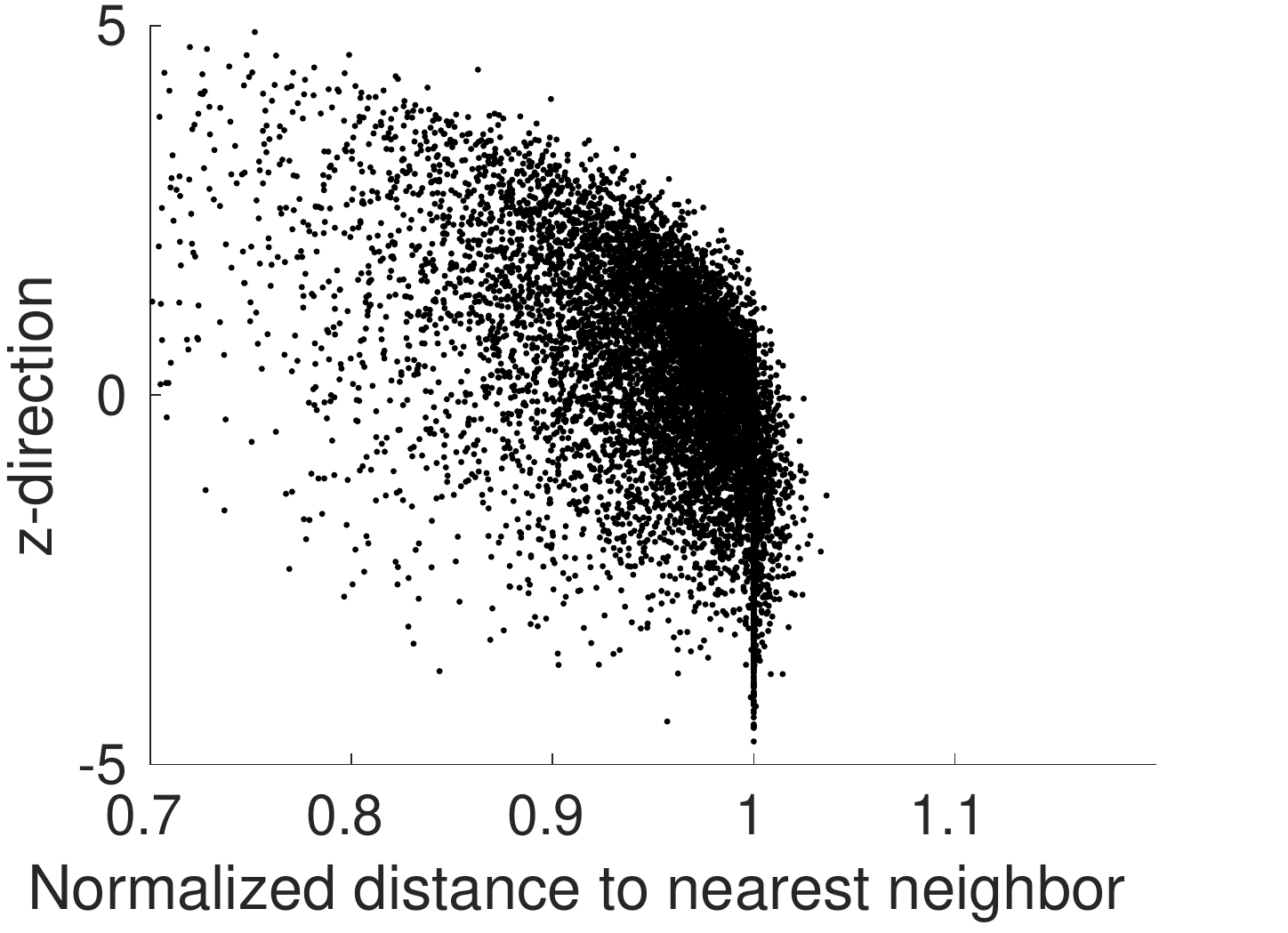}
		\caption{Original}
		\label{fig:zoriginal}
	\end{subfigure}
	\begin{subfigure}{0.325\textwidth}
		\centering
		\includegraphics[width=0.95\textwidth]{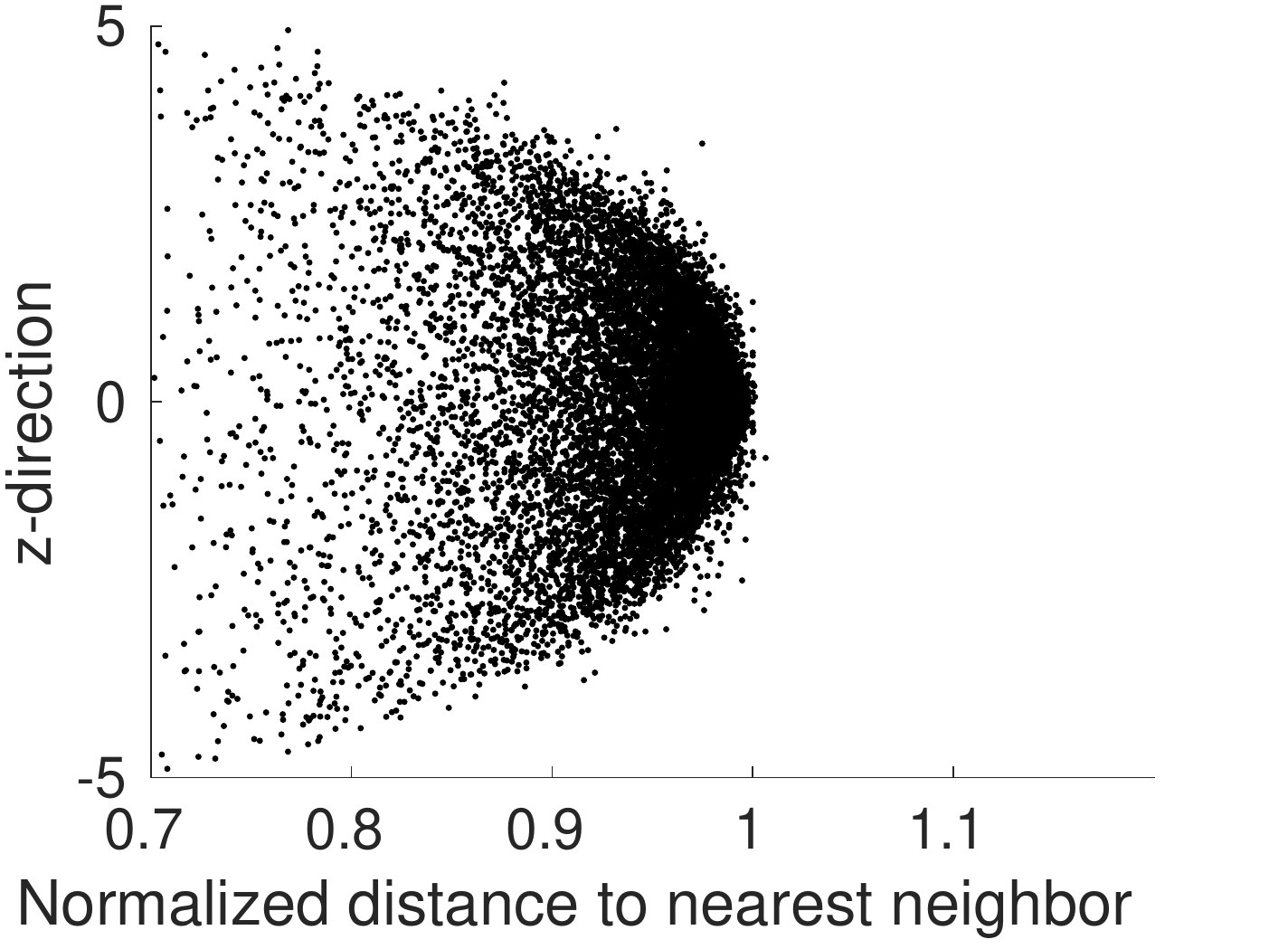}
		\caption{Radially built (1)}
		\label{fig:zcorrected}
	\end{subfigure}
	\begin{subfigure}{0.325\textwidth}
		\centering
		\includegraphics[width=0.95\textwidth]{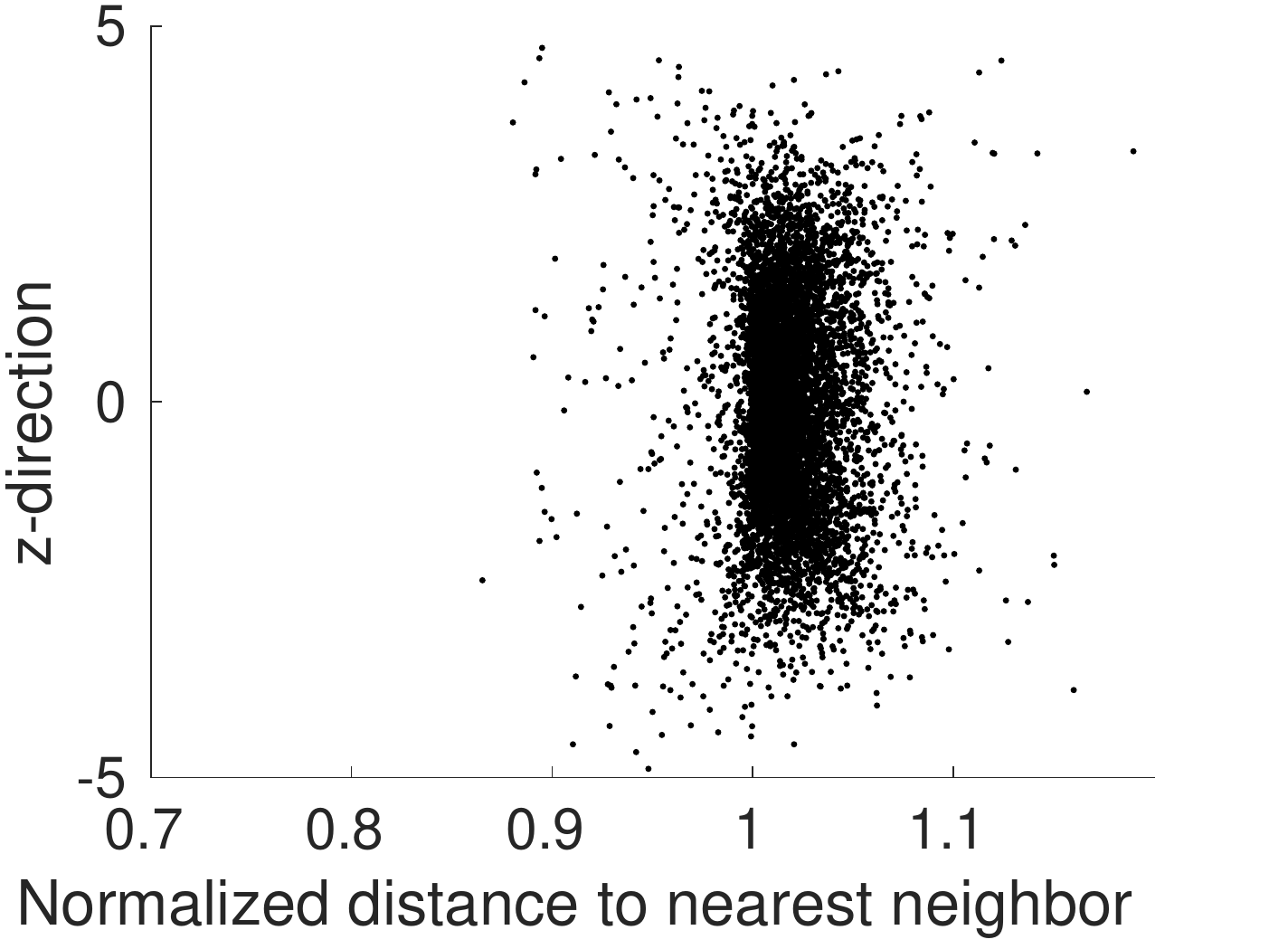}
		\caption{Radially built (2)}
		\label{fig:zradialbuild}
	\end{subfigure}
	\caption{Comparison of the normalized distance to nearest neighbor between the original method and two radially built node set.}
	\label{fig:comparisoninz}
\end{figure}

To investigate the bias in the $z$-direction, the distance to nearest neighbor is plotted as a function of $z$ and compared between node sets in \cref{fig:comparisoninz}. As previously, the distance to nearest neighbor is normalized by dividing by the desired exclusion radius. Here we only compare the original algorithm to both radially built variations. One can see how the radially built node sets avoid the bias seen in the original algorithm and how imposing different minimal spacing requirements through the exclusion radius can affect the distribution of nodes.

\begin{figure}[htbp]
	\centering
	\begin{subfigure}{0.325\textwidth}
		\centering
		\includegraphics[width=0.95\textwidth]{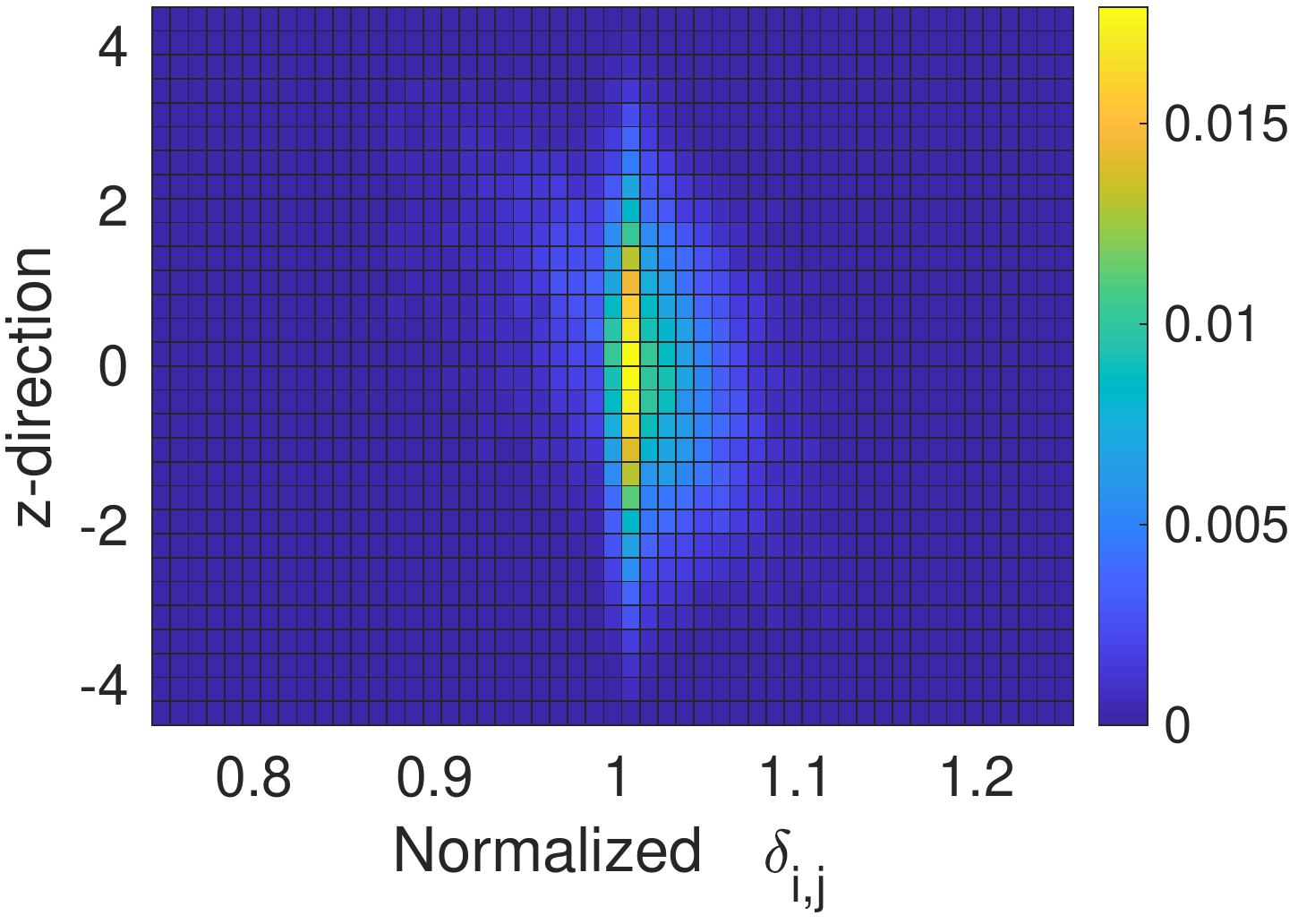}
		\caption{Original}
	\end{subfigure}
	\begin{subfigure}{0.325\textwidth}
		\centering
		\includegraphics[width=0.95\textwidth]{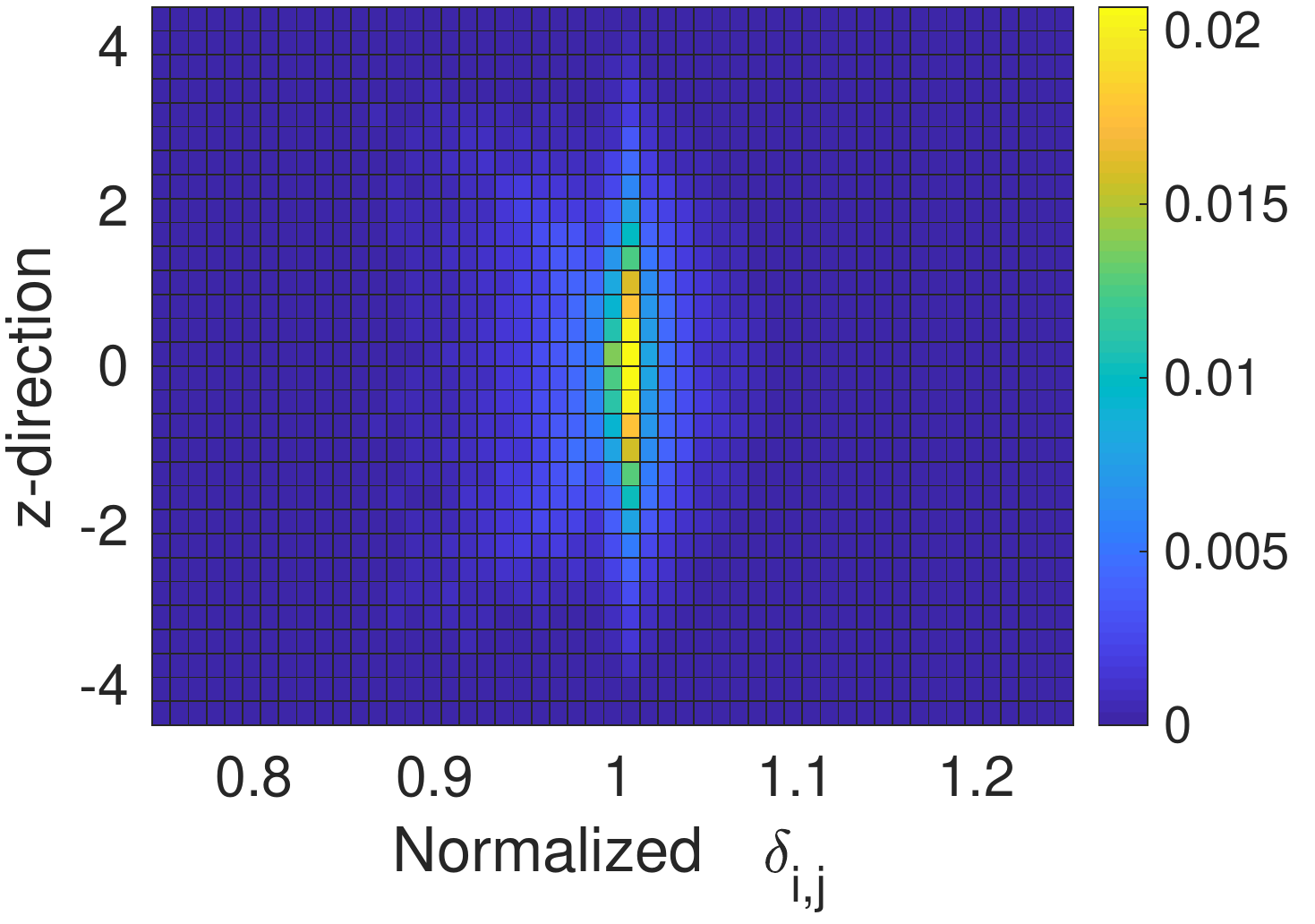}
		\caption{Radially built (1)}
	\end{subfigure}
	\begin{subfigure}{0.325\textwidth}
		\centering
		\includegraphics[width=0.95\textwidth]{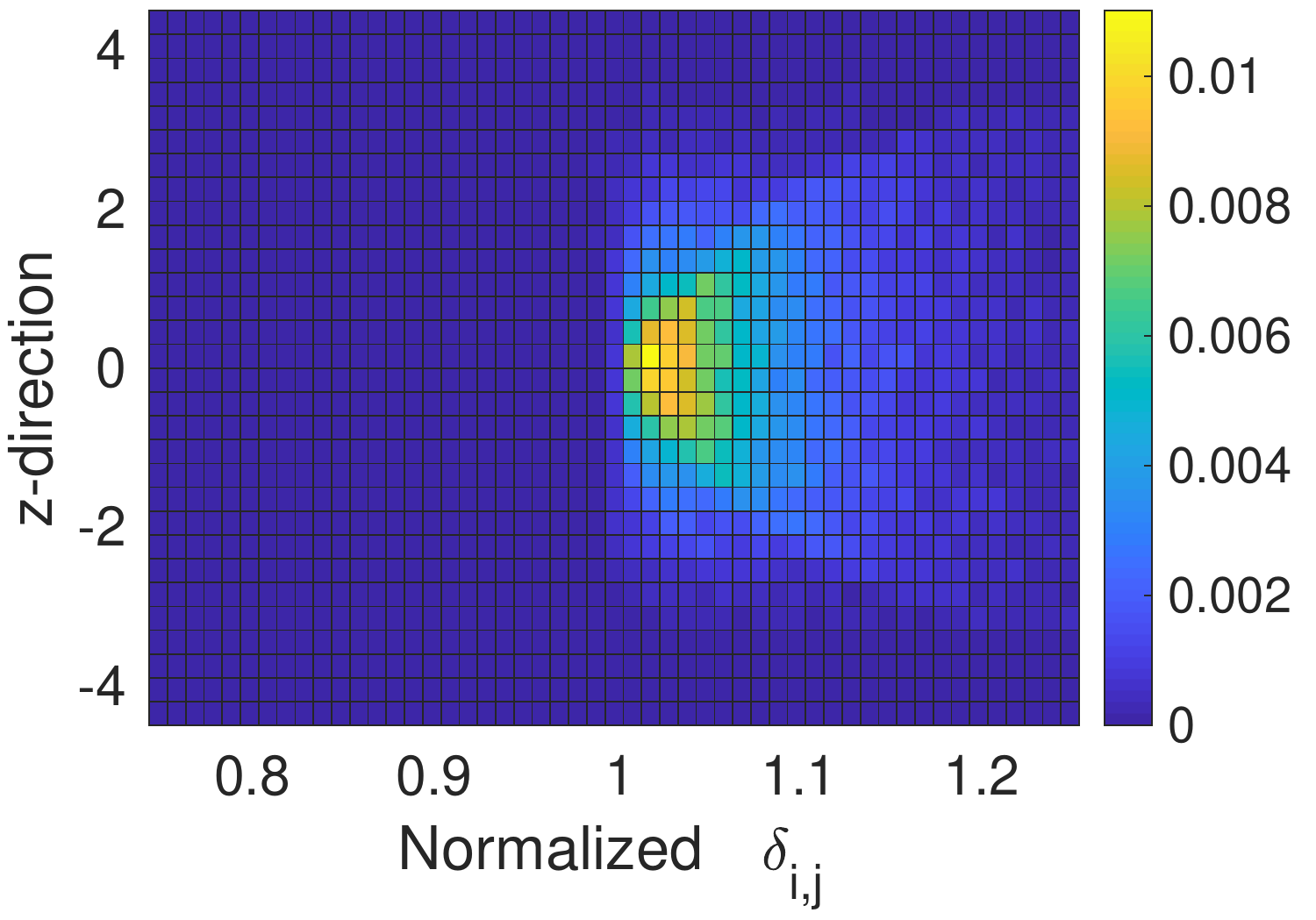}
		\caption{Radially built (2)}
	\end{subfigure}
	\caption{Histogram of normalized distances to the 6 nearest neighbors $\delta_{i,j}$ of each node $p_j$. Color represents the number of counts in each tile normalized by the total number of counts.}
	\label{fig:comparisoninz_6nbrs}
\end{figure}

For additional insight, the distance to $k$ nearest neighbors can be considered as well. In \cref{fig:comparisoninz_6nbrs}, a 2-D histogram shows the normalized distances to the 6 nearest neighbors as a function of $z$. 

\section{Node sets for RBF methods}
\label{sec:RBFmethods}
Here we investigate the application of generated node sets to RBF-FD. The purpose of this work is not to further the development of these meshless methods for solving PDEs. It is instead to present a useful tool for node generation. Hence we will look at the condition number of the collocation matrix as an indicator of the application to RBF-FD methods, and the results of a local interpolation. 

\subsection{Condition Number}

\begin{figure}[htbp]
	\centering
	\begin{subfigure}{0.45\textwidth}
		\centering
		\includegraphics[width=0.9\textwidth]{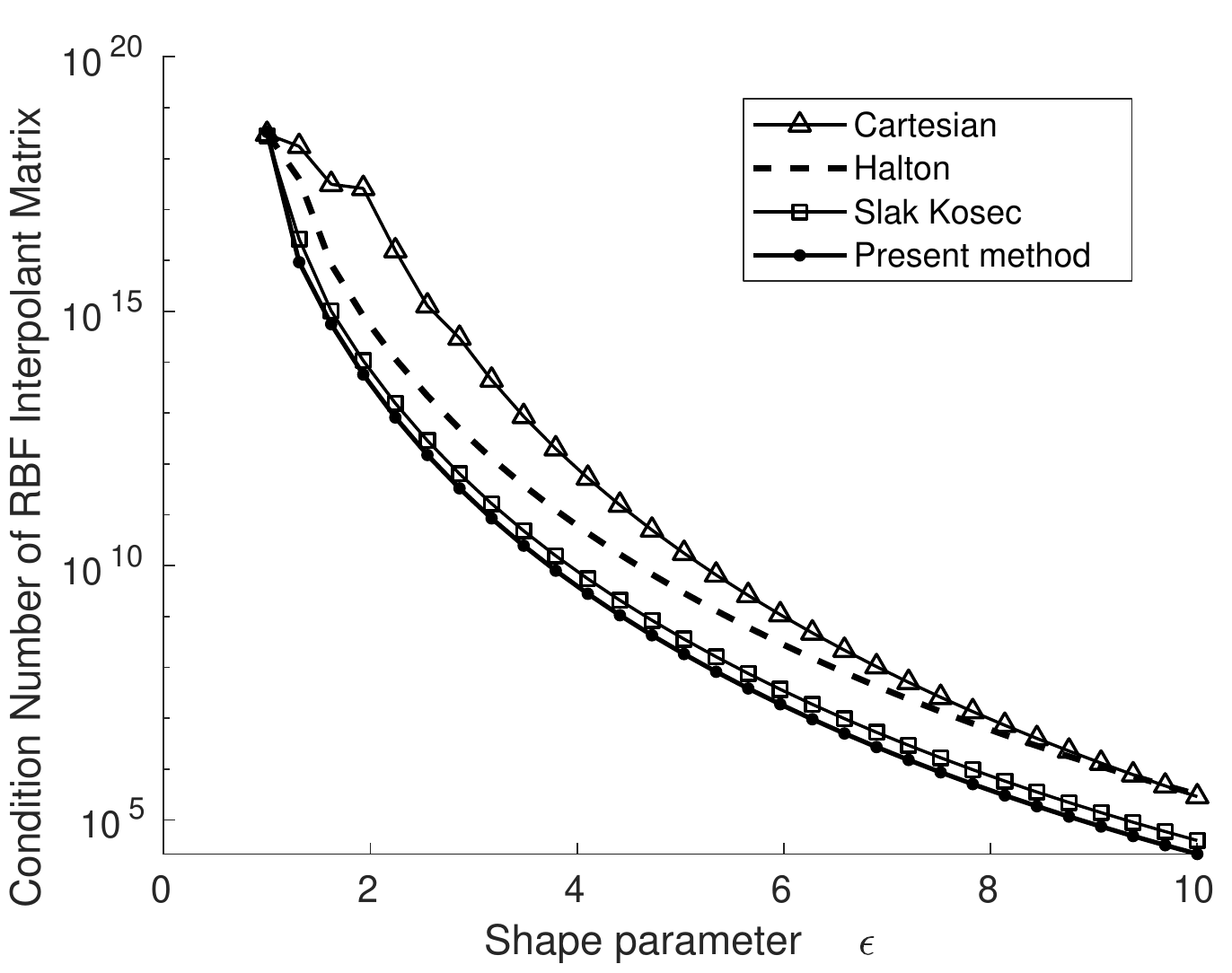}
		\caption{n=80 point stencil and varying $\epsilon$}
		\label{fig:conditionnumber_eps}
	\end{subfigure}
	\begin{subfigure}{0.45\textwidth}
		\centering
		\includegraphics[width=0.9\textwidth]{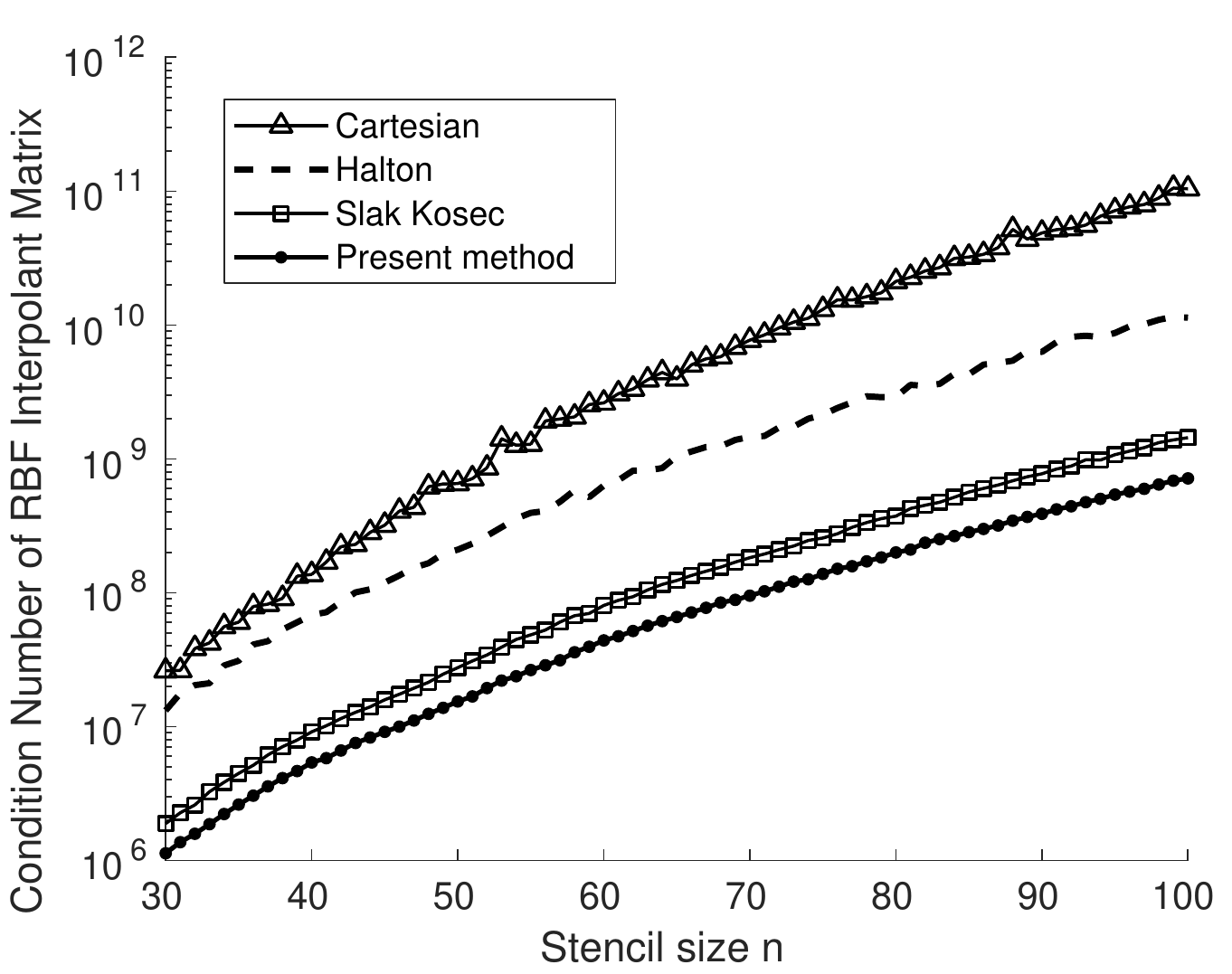}
		\caption{$\epsilon=5$ and varying stencil size $n$}
		\label{fig:conditionnumber_n}
	\end{subfigure}
	\caption{Comparison of the condition number of the RBF-FD matrix using Gaussian RBFs for a 3-D uniform node set of $N \approx 8000$ nodes in the unit cube}
	\label{fig:conditionnumber}
\end{figure}

RBF-FD makes use of a collocation matrix $A$ to obtain the weights for each local stencil of size $n$ \cite{FF2015primer}. One measure of node quality for discretizing PDEs is the condition number of this matrix. The condition number is calculated for a uniform node set of $N \approx 8000$ nodes in the 3-D unit cube. A Gaussian kernel $\phi (r) = e^{-(\epsilon r)^2}$ is used as the basis function, where $r$ is the Euclidean distance from the collocation point and $\epsilon$ is the shape parameter. The condition number is averaged over 300 stencils centered around random points $x_0$ taken from a normal distribution centered in the cube. The result is compared to a Halton node set, a Cartesian lattice, and one generated by the method of Slak \& Kosec in \cref{fig:conditionnumber}. Similar results where the present method and Slak \& Kosec have the lowest condition numbers can be obtained by instead fixing the minimum spacing between node sets (for all except the Halton set) and allowing $N$ to vary. 

It is important to remember that the condition number is not the most useful for measuring node quality, as opposed to the quality metrics investigated in previous sections. In \cite{FF2015primer} it was noted that node irregularity can in some cases reduce condition numbers even while damaging accuracy. When using Gaussian RBFs, a higher condition number can actually give higher accuracy up to a breakpoint where the error spikes. There exist stable algorithms for Gaussian RBFs which bypass these issues in conditioning \cite{FF2015primer,fornberg2013stable}. When using polyharmonic splines augmented with polynomials, as is increasingly popular \cite{Bayona2019role,Bayona2017role}, the condition number becomes irrelevant. 

\subsection{Local Interpolation}
\begin{figure}[htbp]
	\centering
	\begin{subfigure}{0.45\textwidth}
		\centering
		\includegraphics[width=0.9\textwidth]{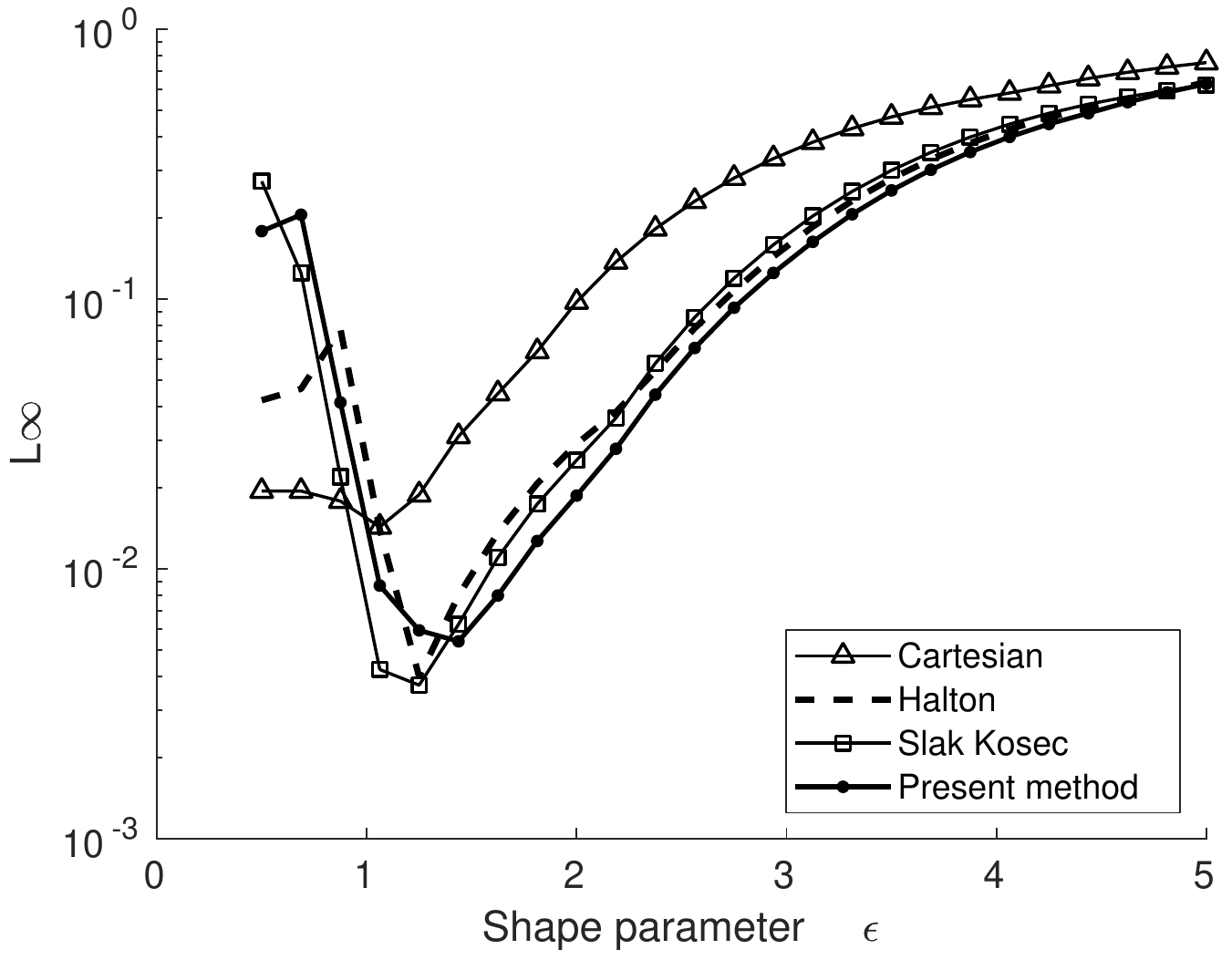}
		\caption{$N \approx 8000$}
	\end{subfigure}
	\begin{subfigure}{0.45\textwidth}
		\centering
		\includegraphics[width=0.9\textwidth]{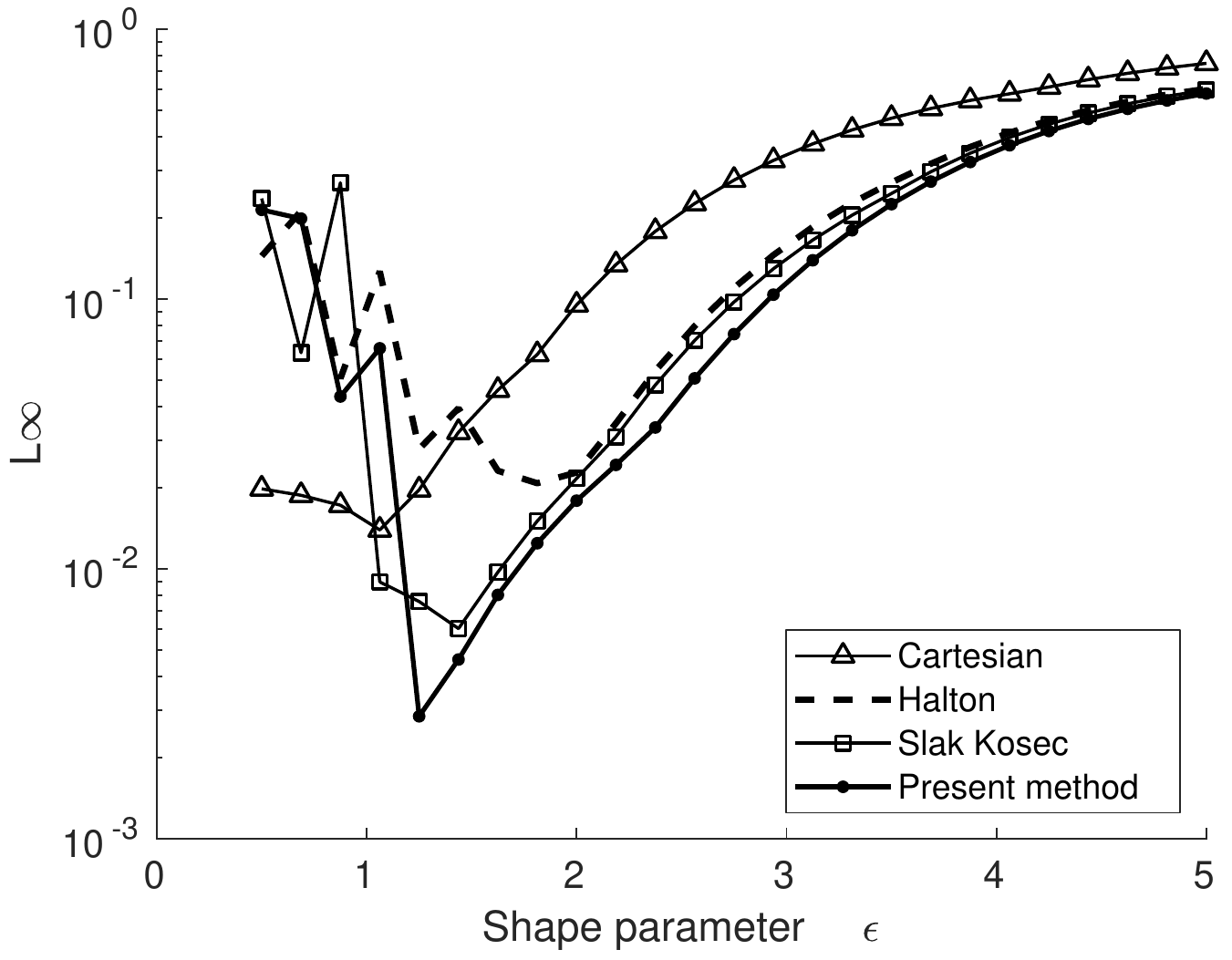}
		\caption{$r = 0.05$}
	\end{subfigure}
	\caption{Error in the local interpolation of $f(R) = \frac{1}{1+R^3}$ where $R$ is distance to the origin.}
	\label{fig:interperror}
\end{figure}

Local interpolation with RBFs provides insight into node quality without getting into the details of solving specific PDEs. We consider a test case of using RBF-FD to calculate a local interpolant to $f(R) = \frac{1}{1+R^3}$ where $R$ is the distance from the origin. Using the same node set in the unit cube, the interpolant was calculated at 100,000 different points using a local stencil of size $n=80$ nodes. The resulting error is compared for different values of the shape parameter $\epsilon$ in \cref{fig:interperror}. The results are shown for both fixed $N \approx 8000$ and fixed minimum spacing $r = 0.05$. For the fixed spacing, we compared to a Halton set with the same number of nodes as the present method. 

\section{Conclusions}
\label{sec:Conclusion}
Methods like RBF-FD for solving PDEs on scattered nodes require that nodes be locally regular and often spatially varying in density. These nodes should satisfy minimum spacing and bounded gaps between nodes. The present method is demonstrated to generate quality node sets in 2-D and 3-D and performs well in comparison to other node sets. It is simple to implement and computationally fast. More complex domains can be treated by generating nodes in a bounded box and then eliminating nodes outside the desired domain. From there, boundary treatment has been discussed in \cite{Fornberg2015nodegen}. Finally, for radially varying density functions a slight modification of the algorithm can allow nodes to be generated in spherical coordinates to reduce bias in the direction of the advancing front. Tests with variable density demonstrated the ability to handle gradients in the radial exclusion function. There may be difficulty in contexts where extreme refinement in small areas of the domain is desired, as this may increase the background grid density to a point where the computational efficiency is lost. 
Future directions include generating nodes from a given boundary set and investigating extensions to adaptive node generation. Source code for the present node generation algorithm in 2-D and 3-D is available at \cite{nodegencode}.

\appendix
\section{Measuring node set quality}
\label{sec:nodesetquality}
There is no general metric of a `good node set', rather, various characteristics are advantageous for different applications. Good point sets for mesh generation or PDE solvers may differ from good points for rendering images or for numerical integration. Low discrepancy is a measure of node quality that has been heavily investigated in relation to numerical integration and Monte Carlo simulations, and has been proposed as a measure of quasi-random node sets \cite{doerr2014calculation}. However, a sequence can have low discrepancy despite having arbitrarily close spacing between nodes: if a pair of points are very close together within a node set of $N$ points, they only add at most $1/N$ to the discrepancy.

Well-spaced nodes have been defined as satisfying a minimal spacing requirement and having bounded gaps \cite{talmor1997well}. These requirements are desirable for mesh-free PDE solvers. The minimal spacing requirement for a uniform density of nodes is clear: given a specified node spacing $r$
\begin{equation}
||p_i - p_j|| \geq r
\end{equation}
for any two distinct nodes $p_i$ and $p_j$.
For the variable density case, the minimal spacing requirement, otherwise known as the empty disk property, is
\begin{equation}
||p_i - p_j|| \geq f(p_i,p_j)
\end{equation}
where $p_i$ is the closest placed node to a new candidate node $p_j$ and $f(p_i,p_j)$ may be one of the variations described by Mitchel et al \cite{mitchell2012variable}:
\begin{equation}
\begin{aligned}
&\text{Prior-disks:} & f(p_i,p_j) &= r(p_i) \\
&\text{Current-disks:} & f(p_i,p_j) &= r(p_j) \\
&\text{Bigger-disks:} & f(p_i,p_j) &= \text{max}(r(p_i),r(p_j)) \\
&\text{Smaller-disks:} & f(p_i,p_j) &= \text{min}(r(p_i),r(p_j)) 
\end{aligned}
\label{eqn:minimalspacingvariations}
\end{equation}

The bounded gaps requirement states that there is an upper bound on the maximum radius of a sphere that can be placed within the node set without including any nodes. As with minimum spacing, this bound should be constant in the case of uniform density nodes, but can be modified for variable density. A node set satisfies the $L$-gap property if for exclusion radius function $r(x)$, the maximum sphere that can be placed within the node set without including any nodes has a radius bounded by $Lr(x)$ where $L$ is a constant \cite{talmor1997well}. 

In the context of RBF-FD it is desirable that nodes be ``locally quasi-uniform", which can be understood intuitively as being roughly equispaced when zoomed in. The term quasi-uniform is well defined on a global sense. A sequence of node sets of size $N$ are globally quasi-uniform if the mesh ratio
\begin{equation}
\gamma_N = \frac{\rho_N}{\delta_N},
\label{eqn:quasiuniformratio}
\end{equation} 
where $\rho_N$ is the covering radius and $\delta_N$ is the maximum distance to the nearest neighbor over the whole set, is bounded as $N \rightarrow \infty$ \cite{hardin2016comparison}. This corresponds to minimizing $\rho$ and maximizing $\delta$ over the whole node set. Although this is a global quality and for variable density node sets one might be interested in looking at the mesh ratio on smaller local patches, it is still always desirable to minimize the global mesh ratio. 

If a Voronoi diagram is constructed from a node set, the covering radius of a node set can be measured as the furthest distance from a node to a vertex of its corresponding Voronoi cell \cite{Conway1993}. Node generation may also be characterized as a sphere packing problem. The sphere packing problem is often separated into a packing problem or a covering problem and a solution to one may not be good for the other. Both can be measured based on a node set's Voronoi diagram. A good packing maximizes the radius of the inscribed circle of the Voronoi cells, while a good covering minimizes the covering radius, which is the radius of the circumscribed circle of the Voronoi cells \cite{Conway1993}. It is natural, therefore to look at the ratio in \eqref{eqn:quasiuniformratio} as a balance between both problems. In using a Voronoi diagram to investigate these metrics, only interior nodes are considered as the Voronoi cells go to infinity at the edges. 

For a node set that satisfies minimal spacing requirements the distance to nearest neighbor is bounded below as $\delta \geq r(x)$. Then the problem of minimizing $\gamma$ can also be reformulated as maximizing the number of nodes in the domain, $N$. To compare further to circle packing, for a uniform node set in 2-D the packing density can be calculated by considering circles around each node, summing the area of the circles within the domain and dividing by the area of the domain. The circles should be half the radius of the exclusion radius. When calculating this packing density, the domain is a box taken from the center of the whole node set in order to avoid boundary effects. It is known that the optimal packing density in the plane is hexagonal, which has a density of $\pi \sqrt{3}/6 \approx 0.9069$. The closer the 2-D packing density is to this, the better. 

A final desirable quality in a node set is local regularity, which requires taking into account the distance to $k$ nearest neighbors. The $k$ neighbors for each node $p_j$ are found and denoted $p_{i,j}$ for $i = 1,2,...k$. The distance to each neighbor is calculated as $\delta_{i,j} = ||p_j - p_{i,j}||$ and an average can be taken over the k neighbors $\bar{\delta_{j}} = \frac{1}{k} \sum_{i=1}^k \delta_{i,j}$ for each node $p_j$. The average $\bar{\delta_j}$ and standard deviation can be taken over the node set as well as the average range of $\max_{j} \delta_{i,j} - \min_{j} \delta_{i,j}$. Again only internal nodes $p_j$ are used to avoid boundary effects.

\bibliographystyle{siamplain}
\bibliography{references}

\end{document}